\documentclass[journal]{IEEEtran}
\usepackage{amsmath}
\usepackage{graphicx}
\usepackage{amssymb}
\usepackage{color}
\usepackage{longtable}
\usepackage{float}
\usepackage{epstopdf}
\usepackage{multirow}
\usepackage{bm}
\usepackage{cite}
\graphicspath{{figure/}}

\newtheorem{defi}{\bf Definition}
\newtheorem{rmk}{\bf Remark}

\newtheorem{prop}{\bf Proposition}
\newtheorem{thm}{\bf Theorem}

\newcommand{\bpara}[4]{ 
\begin{picture}(0,0)%
\setlength{\unitlength}{1pt}%
\put(#1,#2){\rotatebox{#3}{\raisebox{0mm}[0mm][0mm]{%
\makebox[0mm]{$\left.\rule{0mm}{#4pt}\right\}$}}}}%
\end{picture}}

\begin{document}

\title{Minimum Manhattan Distance Approach to Multiple Criteria Decision Making in Multiobjective Optimization Problems}

\author{Wei-Yu~Chiu,~\IEEEmembership{Member,~IEEE}, Gary G. Yen,~\IEEEmembership{Fellow,~IEEE}, and Teng-Kuei Juan  
\thanks{This work was supported by the Ministry of Science and Technology of Taiwan
under Grant 102-2218-E-155-004-MY3. (Corresponding author: W.-Y. Chiu.)}%
\thanks{W.-Y. Chiu and T.-K. Juan are with the Multiobjective Control Lab, Department of Electrical Engineering, Yuan Ze University, Taoyuan 32003, Taiwan (email: chiuweiyu@gmail.com; s1034650@mail.yzu.edu.tw).}%
\thanks{G. G. Yen is with the School of Electrical and Computer Engineering, Oklahoma State University, Stillwater, Oklahoma 74078, USA (email: gyen@okstate.edu).}%
\thanks{
\copyright 2016 IEEE. Personal use of this material is permitted. Permission from IEEE must be
obtained for all other uses, in any current or future media, including
reprinting/republishing this material for advertising or promotional purposes, creating new
collective works, for resale or redistribution to servers or lists, or reuse of any copyrighted
component of this work in other works.
}
\thanks{Digital Object Identifier 10.1109/TEVC.2016.2564158
}
}
\maketitle

\begin{abstract}
A minimum Manhattan distance (MMD) approach to multiple criteria decision making in multiobjective optimization problems (MOPs) is proposed.
The approach selects the finial solution corresponding with a vector that has the MMD from a normalized ideal vector.
This procedure is equivalent to the knee selection described by a divide and conquer  approach that involves iterations of pairwise comparisons.
Being able to systematically assign weighting coefficients to multiple criteria,
the MMD approach is equivalent to a weighted-sum approach.
Because of the equivalence, the MMD approach possesses rich geometric interpretations
that are considered  essential in the field of evolutionary computation.
The MMD approach is elegant because all evaluations can be performed by efficient matrix calculations
without iterations of comparisons.
While the weighted-sum approach may encounter an indeterminate situation in which
a few solutions yield almost the same weighted sum, the MMD approach is able to
 determine the final solution  discriminately.
 Since existing multiobjective evolutionary algorithms
aim for \emph{a posteriori} decision making, i.e., determining the final solution after a set of Pareto optimal solutions is available,
the proposed MMD approach can be combined with them to form a powerful solution method of solving MOPs.
Furthermore, the approach enables scalable definitions of the knee and knee solutions.
\end{abstract}


\section{Introduction}
Multiple Criteria Decision Making (also termed multicriteria decision making, MCDM) or multiple attribute  decision making
occurs naturally in various real-world
problems, e.g., recruitment of employees~\cite{Hamming}, path planning of humanoid robots~\cite{path_follow}, and
factory layout selection for efficient production~\cite{factor_effi}, to name a few.
MCDM is a sub-discipline of operations research~\cite{ORs}.
For an MCDM process, a decision maker (DM) needs to select
a solution (sometimes termed a design, an alternative, or a candidate) out of a set of alternatives based on associated multiple criteria (or attributes). This process can be critical when it involves high stakes, such as a business investment or the sustainability  of a company.

In the field of evolutionary computation, we  encounter MCDM when applying a multiobjective evolutionary algorithm (MOEA) to  solve a multiobjective optimization problem (MOP) or a many-objective optimization problem (MaOP) if more than three objectives are involved.\footnote{In fact, most existing MOEAs,  e.g., PAES~\cite{Knowles99PAES}, PESA~\cite{Corne00PESA}, NSGA-II~\cite{Deb02NSGA-II}, SPEA2~\cite{Zitzler02SPEA2}, and MOEA/D~\cite{Zhang07MOEA}, are designed so that the DM can make a decision after a set of solutions is found.}
By solving an MOP, an approximate Pareto set (APS) and an approximate Pareto front (APF) can be obtained.
The APS consists of  Pareto optimal solutions (or nondominated solutions).
Vectors on the APF correspond with criterion values of Pareto solutions.
A final solution is selected out of the APS based on the performance represented by the APF.
This task could be challenging when the size of the APS  is large.

The selection of a final solution among an APS can be referred to as MCDM in MOPs.
This selection generally depends on the use of a preference model.
In the literature, a preference model can enter the solving process of an MOP at three different stages:
before (\emph{a priori}),  during (progressive), and after (\emph{a posteriori}) the process~\cite{three_category}.
For an \emph{a priori} setting, a series of single-objective formulations combined with preference are often used to convert multiple objectives into one objective~\cite{decom_series}.
For  progressive optimization, the preference of the DM is incorporated into the solution search process~\cite{Rachmawati06,pre_based,branke2004finding,Cvet2002pre,Rach09knee,knee_progress}.
In this case, the size of APS is reduced, leading to a smaller set of candidates, and a further operation is required to determine a final candidate out of the reduced set.
For an \emph{a posteriori} setting, the optimization process is separated from the decision making.

Regarding the construction of a preference model,
one of the most popular ways is to use weighting coefficients or other numerical values that reflect the preference of the DM.
However, related methods can suffer from at least one of the following drawbacks.
They may require careful function normalization and
can be sensitive to the shape of APFs~\cite{marler2004survey,das1997closer}.
They can heavily depend on subjective inputs (or complete knowledge) from the DM~\cite{coit2004system} and, therefore, may devalue information hidden in the APF.
In addition, quantifying the preference of the DM can be difficult.
Even if not impossible, producing preference values can impose much burden on the DM, particularly when
a large number of objectives are involved. Furthermore, some existing approaches lack geometric interpretations (or visualization)
that are considered essential in the evolutionary computation field~\cite{visual_Tusar,geo_im}.\footnote{The importance of geometric interpretations can be readily observed.
For example,  APFs obtained from solving MOPs are often assessed in terms of
the maximum spread, generational distance, and spacing, which have vivid geometric interpretations~\cite{metric}; and
for MaOPs, visualization approaches that address high-dimensional APFs, such as parallel coordinates~\cite{ParCoor,para_cor1,para_cor2}, heatmap~\cite{Pryke15heatmap}, Sammon mapping~\cite{Barton07SammonMapping}, radial coordinate visualization~\cite{Hoffman97RadialCoordinateVisualization}, reduced polar coordinate plot~\cite{He15polar}, self organizing map~\cite{Kohonen01selforganizingmap,SelfOrganizingMaps}, and isomap~\cite{Tenenbaum00isomap},  have become increasingly popular.}

For a preference model that has a geometric interpretation, many studies suggest
the use of knees: solutions corresponding with vectors that geometrically lie in the knee region of the APF should be adopted~\cite{boyd_linear,Mie99knee,Matt04filter,deb2003multi,IET_CTA_14,Rach09knee,das1999characterizing,branke2004finding,Rachmawati06}.
From a geometric perspective, if the shape of the APF is bent, then knee solutions represent those designs that can improve overall performance while sacrificing an insignificant level of performance in certain dimensions~\cite{TSG_15,IET_CTA_14}. In other words, they can exhibit significant improvement in some objectives at the cost of insignificant degradation in the other objectives~\cite{Rach09knee}.

As a preference model, knee selection
has been mostly used before or during the solving process of an MOP.
In~\cite{das1999characterizing}, the knee was associated with the solution to a nonlinear programming problem that was derived from normal-boundary intersection.
It was further shown that the knee is equivalent to a solution of a weighted-sum (WS) problem.
However, this approach requires \emph{a priori} information in practical implementation, and the equivalence was established based on the differentiability of the objective functions, where
 the differentiability  is generally not guaranteed in real-world situations.
In~\cite{branke2004finding}, angle-based and utility-based preference models were proposed and used during the optimization.
While the angle-based model is only suitable for two objectives, the utility-based  model can be extended to any dimensions but
requires a set of weighting parameters. Although the weights can be assigned by sampling, it is not clear how to use the model to uniquely define the knee
in a theoretical framework.
In~\cite{Rach09knee}, gains of improvement and deterioration were evaluated during pairwise comparisons of solutions,
and the knee was characterized as the one that maximized a ratio of improvement over deterioration.
It was argued that such characterization was equivalent to the normal-boundary intersection, but no rigorous proof was provided.

In this study, we primarily focus on the MCDM at the \emph{a posteriori} stage
and propose a minimum Manhattan distance (MMD) approach pertaining to
 knee selection, an appropriate choice of a preference model because of its advantageous properties.
In this case a DM can make a posterior decision and  the resulting MCDM approach can be combined with most MOEAs to form a powerful solution method.
In general, the stage at which a preference model enters the solving process depends on the scenario a DM encounters.
It is not necessary that using a preference model at one stage, e.g., the \emph{a posteriori} stage, is better than that at the other stages, e.g., the \emph{a priori} stage.
Although MOEAs and the MMD approach can be concatenated,
the proposed approach is independent of the choice of MOEAs.
 This is because the decision making process is separated from the optimization process. If two different MOEAs produce the same APF, then the approach yields the same results.
In other words, it is the set of criteria and alternatives that affects the MCDM performance.

The MMD approach determines the final solution associated with the point that has the MMD from an ideal vector.
It can be regarded as a WS approach in which the maximum spread of the APF in each dimension
contributes to weighting coefficients.
In our analysis, we show the equivalence between the WS approach and knee selection described by a divide and conquer (D\&C) approach.
While the MMD approach, WS approach, and  D\&C approach are theoretically equivalent,
the MMD approach is preferred in practice.
In contrast with the D\&C approach that involves iterations of pairwise comparisons for Pareto solutions, the MMD approach can be numerically implemented using efficient matrix calculations.
The WS approach can be affected by the situation in which
 one term in a weighted sum dominates the remaining terms  because of largely distinct scales in objective functions.
 In that case the WS can have difficulty searching for the final solution.
 By contrast, the MMD approach has each term  lie within the interval~$[0,1]$, avoiding this difficulty.

The established equivalence and related analyses enable
the MMD approach to possess the following features:
first, it needs no prior information and avoids using heuristic preference values prescribed by the DM;
second, it has rich geometric interpretations and can be derived from knee selection;
third, it enables a theoretical framework that connects the knee selection with WS approaches;
fourth, it can be analyzed and applied in general situations, which implies that differentiability of objective functions is not required;
 and finally, it allows us to rigorously define the knee and knee solutions, yielding scalable definitions in MaOPs.

The main contributions of this paper are summarized as follows.
We propose an MMD approach to MCDM that has rich geometric interpretations and physical meanings.
We theoretically establish the equivalence between the  MMD approach, WS approach, and  D\&C approach.
Scalable definitions of the knee and knee solutions are rigorously provided.

The rest of this paper is organized as follows. Problem formulations and knee selection are discussed in Section~\ref{sec_moti}.
 The D\&C approach is developed in
Section~\ref{sec_DnC}.
Section~\ref{sec_MMD} presents the MMD approach by connecting it with the WS and D\&C approaches.
Numerical results in Section~\ref{sec_sim} illustrate the effectiveness and efficiency of the MMD approach.
Finally, Section~\ref{sec_con} concludes this paper by addressing the necessity and validity of the proposed methodology.

\section{Problem Formulations and Knee Selection}\label{sec_moti}

In this section, we present mathematical formulations of an MOP and the associated MCDM problem.
We consider knee selection as a way to realize MCDM because of its importance and frequent use in the field of engineering.
Arguments about knee selection from a geometric perspective are offered to facilitate further derivation of the associated algebraic formula.

Consider an MOP
\begin{equation}\label{eq_MOP_ori}
    \begin{split}
      \min & \; \bm{f} (\bm{x})\\
      \mbox{ s.t. }  & \bm{x} \in \Omega
    \end{split}
\end{equation}
where
\begin{equation*}
\bm{f} (\bm{x})=
\left[
 \begin{array}{cccc}
 f_1(\bm{x}) & f_2(\bm{x}) & ... & f_N(\bm{x}) \\
 \end{array}
\right]^T
\end{equation*}
represents  a vector of objective functions, $\bm{x}$ represents a vector of decision variables, and
$\Omega$ denotes the feasible search space. The optimality in~(\ref{eq_MOP_ori}) is often defined by Pareto dominance~\cite{MOEA_bk1,MOEA_bk2}.

\begin{defi}[Pareto dominance]\label{def_domi}
In the decision variable space of~(\ref{eq_MOP_ori}), a point
$\bm{x}' \in \Omega$ dominates another point $\bm{x}'' \in \Omega$
if  the conditions  $f_i(\bm{x}')  \leq f_i(\bm{x}''),i=1,2,...,N,$
hold true and at least one inequality is strict.
 In this case, we denote  $\bm{x}' \preceq   \bm{x}'' $.
A point that is not dominated by other points is termed a nondominated point.
\end{defi}

\begin{defi}[Pareto optimal set]\label{def_P_opt}
 The Pareto optimal set $\mathcal{PS}$ of~(\ref{eq_MOP_ori}) is defined as the set of all nondominated points, i.e.,
 \begin{equation*}
  \mathcal{PS}=  \{  \bm{x} \in \Omega : \nexists \bm{x}' \in  \Omega  \mbox{ s.t. } \bm{x}' \preceq
   \bm{x} \}.
 \end{equation*}
\end{defi}

\begin{defi}[Pareto front]\label{def_PF}
The Pareto front (PF) of~(\ref{eq_MOP_ori}) is defined as the image of the Pareto optimal set based on  the mapping of the vector-valued objective function $\bm{f}$, i.e.,
the set $\bm{f}(\mathcal{PS})$ is the PF.
\end{defi}

Without loss of generality, in~(\ref{eq_MOP_ori})
only minimization  is considered  because maximizing an objective function can be equivalently transformed into minimizing the negative value of the objective function.
An MCDM problem occurs when we apply an MOEA to solve~(\ref{eq_MOP_ori}).
After the solving process, we can obtain an APS ($\mathcal{PS}^{A}$) and corresponding APF ($\mathcal{PF}^{A}$), denoted by
\begin{equation}\label{eq_APF_APS}
\begin{split}
  \mathcal{PS}^{A}{ }={ }  & \{ \bm{x}_1, \bm{x}_2, ... , \bm{x}_M \}  \mbox{ and }\\
 \mathcal{PF}^{A}{ }={ }  &  \{ \bm{f}(\bm{x}_1) , \bm{f}(\bm{x}_2), ... , \bm{f}(\bm{x}_M)\}
\end{split}
\end{equation}
respectively. The MCDM problem is about how to select one solution from $\mathcal{PS}^{A}$
based on information hidden in $\mathcal{PF}^{A}$.
In practice, $N$ or $M$ can be large and
thus MCDM can be a challenging task.

In this study, we are interested in developing a method for determining a final solution that corresponds with a vector in the knee region of the APF.
Such a method can possess a geometric interpretation and avoid heuristic assignment of weighting coefficients.
To begin with, we consider a simple scenario in which only two objectives are involved, i.e., $N=2$.
Since finding a way to compare two solutions should be easier than
developing a method for comparing all solutions simultaneously,
the MCDM problem is divided into several subproblems, each of which
consists of only two solutions from $\mathcal{PS}^{A}$.

\begin{figure}
 \centering
\begin{equation*}
 \begin{array}{cc}
 \includegraphics[width=5cm]{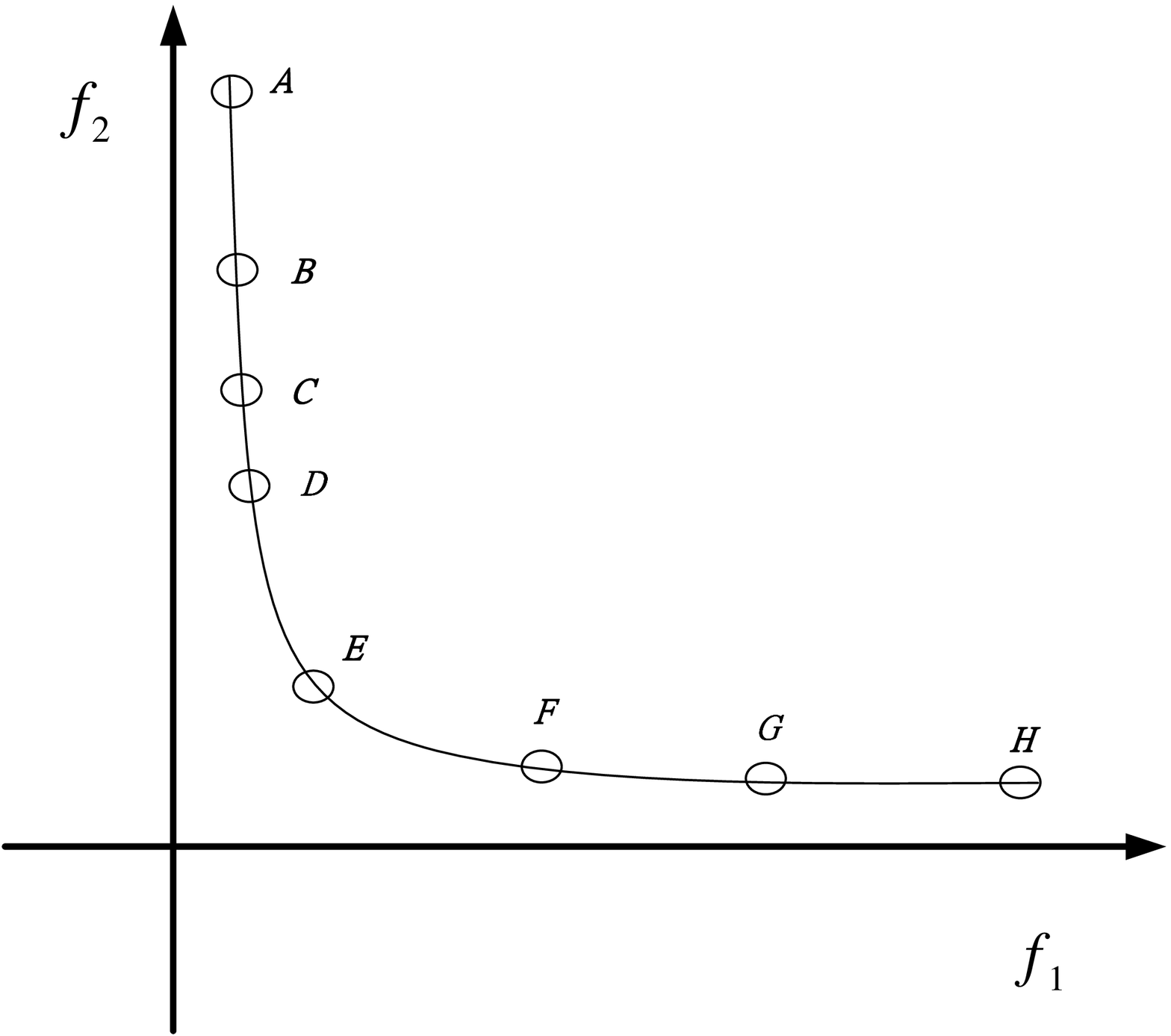} &   \includegraphics[width=3.6cm]{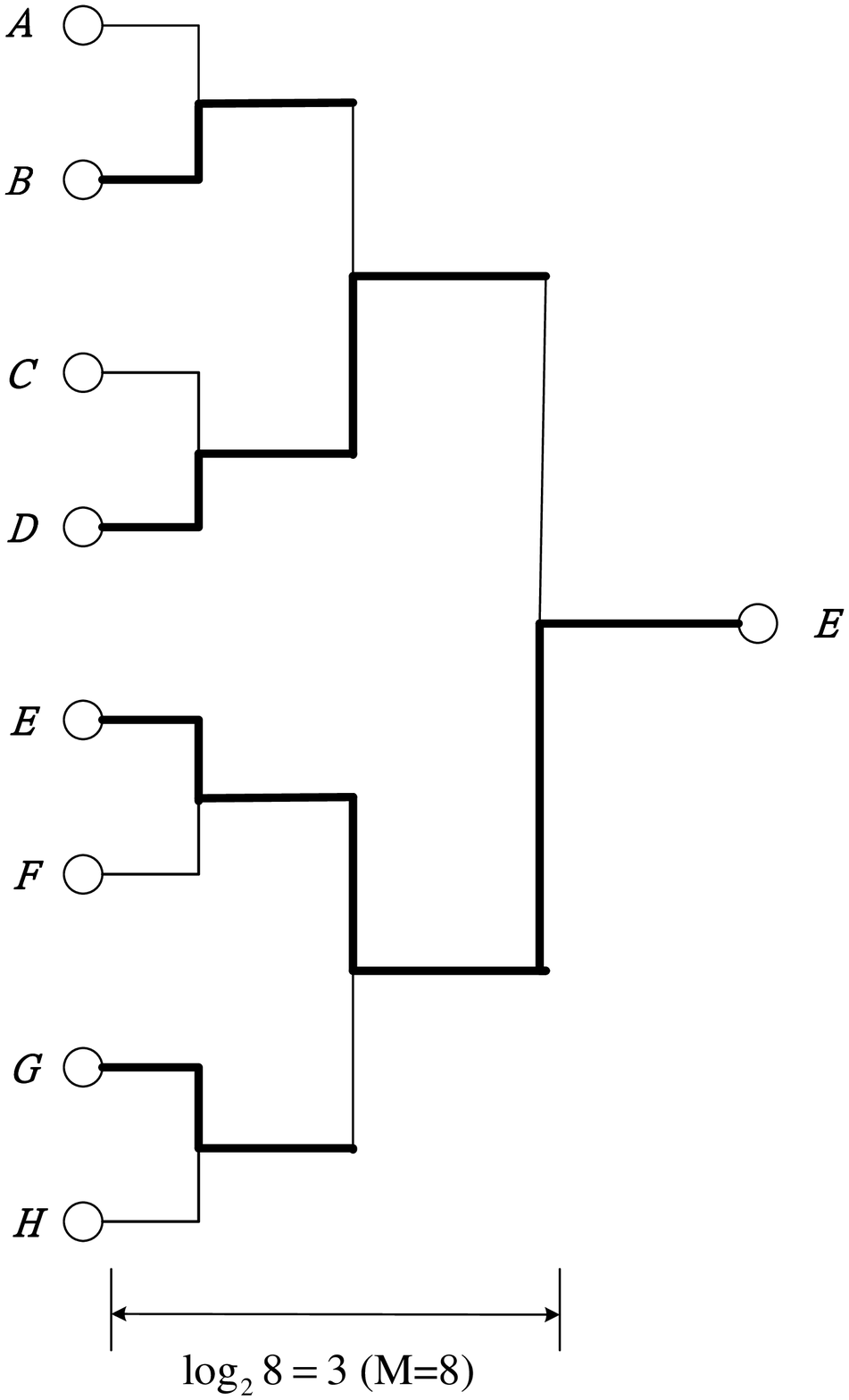}\\
   \mbox{(a)} & \mbox{(b)} \\
 \end{array}
\end{equation*}
  \caption{Demonstration of pairwise comparison between vectors on an APF. (a) Points $\bm{A}, \bm{B}, ..., \bm{G}$, and $\bm{H}$ are nondominated.
   (b) For $M=8$ solutions, $\log_2 8$ iterations of comparisons are performed to obtain a final solution.
  }\label{fig_PF_point}
\end{figure}

In Fig.~\ref{fig_PF_point}(a), the PF is symbolically represented by the solid curve.
Points $\bm{A}$, $\bm{B}$, ..., $\bm{G}$ and $\bm{H}$ are vectors marked on the $\mathcal{PF}^{A}$.
Comparing points~$\bm{A}$ with $\bm{B}$, we note that point~$\bm{A}$ gives a better $f_1$ value but a worse  $f_2$ value.
While the difference between points $\bm{A}$ and $\bm{B}$ in terms of the $f_1$ value is insignificant,
point~$\bm{B}$ has a substantially better $f_2$ value than point~$\bm{A}$.
It is thus reasonable to select point $\bm{B}$  when points $\bm{A}$ and $\bm{B}$ are compared.
This selection can be interpreted as follows:
transition from point~$\bm{A}$ to point~$\bm{B}$ is preferred because the substantial percentage improvement in the $f_2$ dimension can outweigh the insignificant performance degradation in the $f_1$ dimension.
Therefore, when an APS is available, we can divide solutions into pairs for comparison in the objective function space. The number of candidate solutions
can be shrunk by half after each iteration of pairwise comparisons. Ideally, the final solution can be obtained after
$\log_2 M$ iterations, as illustrated in Fig.~\ref{fig_PF_point}(b).

Geometrically, a neighborhood of point~$\bm{E}$ can constitute a region termed a knee region.
Solutions corresponding with vectors in the knee region are termed knee solutions. These solutions are of interest
because they achieve an excellent level of overall performance while sacrificing each objective to a small extent.
The way we chose point~$\bm{E}$ follows the idea of knee selection.
To extend our method to a higher dimension, i.e., a larger $N$, we generalize our selection philosophy from the case where $N=2$:
for a pairwise comparison, point~$\bm{B}$ is preferred to point~$\bm{A}$ if transition from~$\bm{A}$ to~$\bm{B}$ yields a larger improvement percentage  in one dimension than the degradation percentage  in the other dimension.
For a pairwise comparison when $N\geq 3$, we argue that point~$\bm{B}$ is preferred to point~$\bm{A}$ if
transition from $\bm{A}$ to $\bm{B}$ yields a positive net improvement percentage.

To illustrate our selection philosophy, we consider the case in which $N=3$.
Suppose that an APF is available and that
transition from point~$\bm{A}$ to point B yields a triple $(-20\%,15\%,15\%)$, where the negative sign represents
the performance degradation.
In this example, we prefer $\bm{B}$ to $\bm{A}$ because the net improvement percentage  is $-20\%+15\%+15\%= 10\%$, which is positive.

We have derived a selection strategy from pairwise comparisons.
Solutions are divided into groups, each of which consists of only two elements.
For any pairwise comparison within a group, one is preferred if the net improvement percentage  is positive.
This approach is termed the  D\&C and considered as knee selection. There are two features  the approach must possess to ensure a legitimate MCDM.
First,  the transition direction should not affect the selection, which implies that
the net improvement percentage from point $\bm{A}$ to $\bm{B}$ should  equal
the negative value of the net improvement percentage from point~$\bm{B}$ to~$\bm{A}$.
Otherwise, we can encounter a dilemma in which one solution is preferred in one transition direction but
 not preferred in the other transition direction. Second, the D\&C approach should be able to address  the situation in which
the net improvement percentage equals zero. In this case, there is no reason to move from one point to the other, leading to incomparability.
To formalize the  D\&C approach that possesses these two features, we need rigorous mathematical definitions
and a preference model.

\section{Divide and Conquer Approach}\label{sec_DnC}

In this section, we  define mathematically the net improvement percentage, leading to
 a  preference model for pairwise comparisons.
By using this model, two solutions can be compared even if their associated vectors on the APF are close to each other.
The concept of equivalence class is then introduced to address the situation in which the net improvement percentage equals zero.
Finally, we show that the D\&C approach can produce a unique and consistent class. Therefore, the final result is independent of how pairwise comparisons are conducted.
The existence of the unique class enables us to  rigorously define the knee in the APF and knee solution in the APS.

Referring to the notations in~(\ref{eq_MOP_ori}) and (\ref{eq_APF_APS}), we define the improvement percentage  and net improvement percentage  as follows.

\begin{defi}\label{def_IP}
For a transition from solution $\bm{x}_i$ to solution $\bm{x}_j$, denoted by $\bm{x}_i \rightarrow \bm{x}_j$,  the improvement percentage  in the $n$th dimension, denoted by $IP_n(\bm{x}_i \rightarrow \bm{x}_j)$,  is defined as
\begin{equation}\label{eq_PI}
  IP_n(\bm{x}_i \rightarrow \bm{x}_j) =  \frac{f_n(\bm{x}_i) - f_n(\bm{x}_j)    }{L_n}  \times 100 \%
\end{equation}
where
\begin{equation*}
  L_n= \max_m f_n(\bm{x}_m) -  \min_{m} f_n(\bm{x}_{m}).
\end{equation*}
\end{defi}

\begin{defi}\label{def_NIP}
For a transition from solution $\bm{x}_i$ to solution $\bm{x}_j$, $IP(\bm{x}_i \rightarrow \bm{x}_j)$ denotes the net improvement percentage   and is defined as
\begin{equation}\label{eq_PI}
  IP(\bm{x}_i \rightarrow \bm{x}_j) = \sum_{n=1}^N   IP_n(\bm{x}_i \rightarrow \bm{x}_j).
\end{equation}
\end{defi}

Definitions~\ref{def_IP} and~\ref{def_NIP} can be used to construct a preference model
for pairwise comparisons:
\begin{equation}\label{eq_pref1}
  \bm{x}_j \mbox{ is preferred to } \bm{x}_i \mbox{ if } IP(\bm{x}_i \rightarrow \bm{x}_j)>0.
\end{equation}
In other words, solution $\bm{x}_j$ is selected instead of solution $\bm{x}_i$ if the net improvement percentage associated with $\bm{x}_i \rightarrow \bm{x}_j$ is positive.
The following theorem shows that  the transition direction does not affect the selection.

\begin{thm}\label{thm_tran_dir}
By using the preference model in~(\ref{eq_pref1}),
 $\bm{x}_j$ is preferred to  $\bm{x}_i$ if and only if
\begin{equation*}
 IP(\bm{x}_j \rightarrow \bm{x}_i)< 0.
\end{equation*}
\end{thm}

{\emph{Proof:}
This can be verified by noting that
\begin{equation}\label{eq_trans_dir}
\begin{split}
  & IP(\bm{x}_i \rightarrow \bm{x}_j)=\sum_{n=1}^N   IP_n(\bm{x}_i \rightarrow \bm{x}_j) \\
  { }={ } &   - \sum_{n=1}^N  IP_n(\bm{x}_j \rightarrow \bm{x}_i)=  -  IP(\bm{x}_j \rightarrow \bm{x}_i).
\end{split}
\end{equation}
Therefore, $IP(\bm{x}_i \rightarrow \bm{x}_j)>0$ if and only if $ IP(\bm{x}_j \rightarrow \bm{x}_i)<0$.
\hfill $\Box$

Although the preference model is mostly valid,   we might encounter a situation in which two solutions $\bm{x}_j$ and $\bm{x}_i$ are incomparable, i.e., $IP(\bm{x}_i \rightarrow \bm{x}_j)=0$. In such a case, we cannot say that one solution is preferred to the other. To avoid possible incomparability, we employ the concept of an equivalence relation to classify two solutions that yield $IP(\bm{x}_i \rightarrow \bm{x}_j)=0$ into a same class~\cite{set_theory}.

\begin{defi}
A relation in $\mathcal{PS}^{A}$  is a subset of $\mathcal{PS}^{A}\times \mathcal{PS}^{A}$, where ``$\times$'' represents the Cartesian product.
Let $G$ denote a relation in  $\mathcal{PS}^{A}$.
$G$ is reflexive if $(\bm{x},\bm{x})\in G$ for all $\bm{x} \in \mathcal{PS}^{A}$;
 $G$ is symmetric if $(\bm{x}_j, \bm{x}_i)\in G  $ implies  $(\bm{x}_i, \bm{x}_j)\in G$; and $G$ is transitive if $(\bm{x}_i,\bm{x}_j) \in G$ and $(\bm{x}_j,\bm{x}_k) \in G$ imply
  $(\bm{x}_i,\bm{x}_k) \in G$.
\end{defi}

\begin{defi}
A relation $G$ in $\mathcal{PS}^{A}$ is an equivalence relation  if
 it is reflexive, symmetric, and transitive. For an equivalence relation $G$,
 we use $\bm{x}_i\simeq \bm{x}_j$ to represent  $(\bm{x}_i,\bm{x}_j)\in G  $.
\end{defi}

By using the definitions of relation, we have the following result.

\begin{thm}\label{thm_euqi_rel}
  Define a relation  in $\mathcal{PS}^{A}$ as follows:
 $\bm{x}_i \simeq \bm{x}_j$ if
\begin{equation}\label{eq_defi_relation}
IP(\bm{x}_i \rightarrow \bm{x}_j)=0.
\end{equation}
 Then the relation represented by $\simeq$ is an equivalence relation.
\end{thm}

{\emph{Proof:}
The relation is reflexive because
\begin{equation*}
\begin{split}
      & IP(\bm{x}_i \rightarrow \bm{x}_i)= \sum_{n=1}^N   IP_n(\bm{x}_i \rightarrow \bm{x}_i)\\
 ={ } & \sum_{n=1}^N \frac{f_n(\bm{x}_i) - f_n(\bm{x}_i)    }{L_n}=0.
 \end{split}
\end{equation*}
The relation  is symmetric because
\begin{equation*}
\begin{split}
        & \bm{x}_i \simeq \bm{x}_j  \\
 \Rightarrow{ }   & IP(\bm{x}_i \rightarrow \bm{x}_j)=0 \\
 \Rightarrow{ }   &   \sum_{n=1}^N \frac{f_n(\bm{x}_i) - f_n(\bm{x}_j)    }{L_n}=0 \\
 \Rightarrow{ }   &   \sum_{n=1}^N \frac{f_n(\bm{x}_j) - f_n(\bm{x}_i)    }{L_n}=0 \\
 \Rightarrow{ }   & IP(\bm{x}_j \rightarrow \bm{x}_i)=0 \\
 \Rightarrow{ }   & \bm{x}_j \simeq \bm{x}_i.
\end{split}
\end{equation*}
Finally, the relation is transitive because
\begin{equation*}
\begin{split}
        & \bm{x}_i \simeq \bm{x}_j \mbox{ and } \bm{x}_j \simeq \bm{x}_k \\
 \Rightarrow{ }   & IP(\bm{x}_i \rightarrow \bm{x}_j)=0 \mbox{ and } IP(\bm{x}_j \rightarrow \bm{x}_k)=0\\
 \Rightarrow{ }   &   \sum_{n=1}^N \frac{f_n(\bm{x}_i) - f_n(\bm{x}_j)    }{L_n}=0 \mbox{ and } \sum_{n=1}^N \frac{f_n(\bm{x}_j) - f_n(\bm{x}_k)    }{L_n}=0 \\
 \Rightarrow{ }   &    \sum_{n=1}^N \frac{f_n(\bm{x}_i) - f_n(\bm{x}_j)    }{L_n} + \sum_{n'=1}^N \frac{f_{n'}(\bm{x}_j) - f_{n'}(\bm{x}_k)    }{L_{n'}}=0 \\
 \Rightarrow{ }   & \sum_{n=1}^N \frac{f_n(\bm{x}_i) - f_n(\bm{x}_k)    }{L_n} \\
 \Rightarrow{ }   & IP(\bm{x}_i \rightarrow \bm{x}_k)=0 \\
 \Rightarrow{ }  & \bm{x}_i \simeq \bm{x}_k.
\end{split}
\end{equation*}
Since the relation is reflexive, symmetric, and transitive, it is an equivalence relation, which completes the proof.\hfill $\Box$

The preference model  in~(\ref{eq_pref1}) does not address the case in which $IP$ is equal to zero.
When the case occurs, the solutions in comparison are considered equivalent and, hence,
we classify these solutions into equivalence classes.

\begin{defi}\label{def_class}
  For an equivalence relation $\simeq$ in $\mathcal{PS}^{A}$, the equivalence class of $\bm{x}_i$, denoted by $[\bm{x}_i]$, is the set
  \begin{equation}\label{eq_equi_class}
    [\bm{x}_i] =\{\bm{x}\in \mathcal{PS}^{A}: \bm{x}_i \simeq \bm{x}   \}.
  \end{equation}
\end{defi}

A well-known result regarding an equivalence relation in a set is as follows~\cite{algebra_set,set_theory}.
\begin{prop}\label{thm_partition}
  Given an equivalence relation $\simeq$ in $\mathcal{PS}^{A}$, equivalence classes induced by $\simeq$ give a partition of $\mathcal{PS}^{A}$.
\end{prop}

By using the equivalence relation defined in Theorem~\ref{thm_euqi_rel}, we can classify a pair of solutions that yield
$IP(\bm{x}_i \rightarrow \bm{x}_j)=0$ into the same equivalence class $[\bm{x}_i]$. A class $[\bm{x}]$ instead of a point $\bm{x} \in \mathcal{PS}^{A}$   is considered  a single mathematical object (or element) afterwards. Since equivalence classes give a partition according to Proposition~\ref{thm_partition}, all solutions in $\mathcal{PS}^{A}$ belong to certain classes.
 To avoid the incomparability problem,
the preference model in~(\ref{eq_pref1}) can be modified as follows:
\begin{equation}\label{eq_pref2}
  [\bm{x}_j] \mbox{ is preferred to } [\bm{x}_i] \mbox{ if } IP(\bm{x}_i \rightarrow \bm{x}_j)>0.
\end{equation}
The following theorem shows that
the most preferred element exists and hence,
the knee in the APF and the knee solution in the APS can be defined accordingly.

\begin{thm}\label{thm_minimal_element}
Let
\begin{equation}\label{eq_set_class}
  \mathcal{PS}^{A}_c=\{ [\bm{x}]: \bm{x} \in  \mathcal{PS}^{A}  \}
\end{equation}
denote the partition induced by the equivalence relation defined in Theorem~\ref{thm_euqi_rel},
and  $[\bm{x}_j]\prec_{\bm{k}} [\bm{x}_i] $ denote that $[\bm{x}_j]$  is preferred to  $[\bm{x}_i]$ based on keen selection.
 There exists an unique element $[\bm{x}^*]\in \mathcal{PS}^{A}_{c} $ such that  $ [\bm{x}^*] \prec_{\bm{k}} [\bm{z}]$ for all $[\bm{z}] \in \mathcal{PS}^{A}_{c} \setminus  \{  [\bm{x}^*] \}$.
\end{thm}

{\emph{Proof:}
Given two distinct $[\bm{x}_i]$ and $[\bm{x}_j]\in \mathcal{PS}^{A}_{c}$, we have either
\begin{equation*}
IP(\bm{x}_i \rightarrow \bm{x}_j)>0 \mbox{ or } IP(\bm{x}_j \rightarrow \bm{x}_i)>0.
\end{equation*}
Therefore, each element in $\mathcal{PS}^{A}_{c}$ is comparable.
Since  $\mathcal{PS}^{A}_{c}$ has finite elements, the existence of $[\bm{x}^*]$ holds true.
The uniqueness is verified by noting that we cannot have
\begin{equation*}
  [\bm{x}^*] \prec_{\bm{k}} [\bm{z}]  \mbox{ and }  [\bm{z}]  \prec_{\bm{k}} [\bm{x}^*]
\end{equation*}
for some  $[\bm{z}]\not=[\bm{x}^*]$
 because conditions
\begin{equation*}
  IP(\bm{x}^* \rightarrow \bm{z})>0  \mbox{ and }  IP(\bm{x}^* \rightarrow \bm{z})< 0
\end{equation*}
 cannot hold true simultaneously.
\hfill $\Box$

Based on Theorem~\ref{thm_minimal_element}, we summarize the D\&C approach as follows:
\begin{enumerate}
  \item Input $\mathcal{PS}^{A}$  and $\mathcal{PF}^{A}$.
  \item Construct $\mathcal{PS}^{A}_c$.
  \item Perform pairwise comparisons among $[\bm{x}]\in \mathcal{PS}^{A}_c$ according to~(\ref{eq_pref2}).
  \item Output the element $[\bm{x}^*]$.
\end{enumerate}
In addition,
because of the uniqueness of $[\bm{x}^*]$, we are able to define the knee and knee solution as follows.
\begin{defi}
  Given $\mathcal{PS}^{A}$  and $\mathcal{PF}^{A}$, the knee solution is the unique element $[\bm{x}^*]$ in Theorem~\ref{thm_minimal_element} and the knee is the set
\begin{equation*}
   \bm{f}([\bm{x}^*]) =\{  \bm{f}(\bm{x}): \bm{x} \in  [\bm{x}^*]   \} \subset \mathcal{PF}^{A}.
\end{equation*}
\end{defi}

This section presented knee selection described by the D\&C approach.
The approach is efficient because the size of $\mathcal{PS}^{A}$ can be reduced by half  after each iteration of pairwise comparisons.
In practice, iterations of comparisons can be replaced by elegant matrix calculations.
We show this by connecting the knee selection with the MMD in the next section.

\section{Minimum Manhattan Distance Approach}\label{sec_MMD}

This section develops the proposed MMD approach to MCDM:
the solution that minimizes the distance from a normalized ideal vector is selected.
The MMD approach originates from the D\&C approach.
First, the D\&C approach is  transformed into a WS approach by rearranging terms in an inequality that is associated with the preference model.
Next, the WS approach is transformed into the MMD approach by adding entries of an ideal vector to the weighted sum.
The established equivalence between these approaches allows the MMD approach to possess rich geometric interpretations.

\begin{figure*}
  \centering
\begin{equation*}
    \begin{array}{cc}
     \includegraphics[width=9cm]{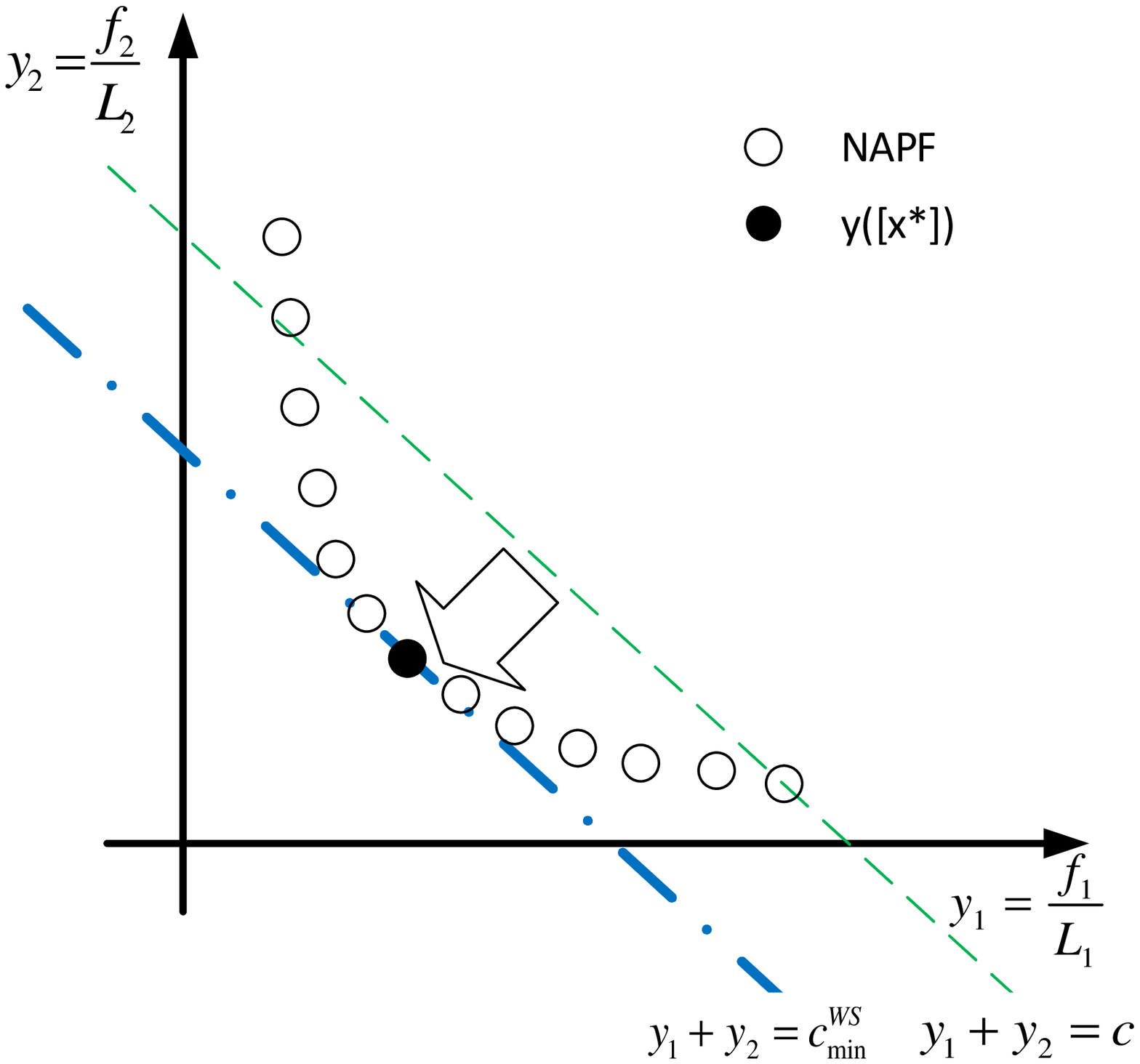} &  \includegraphics[width=9cm]{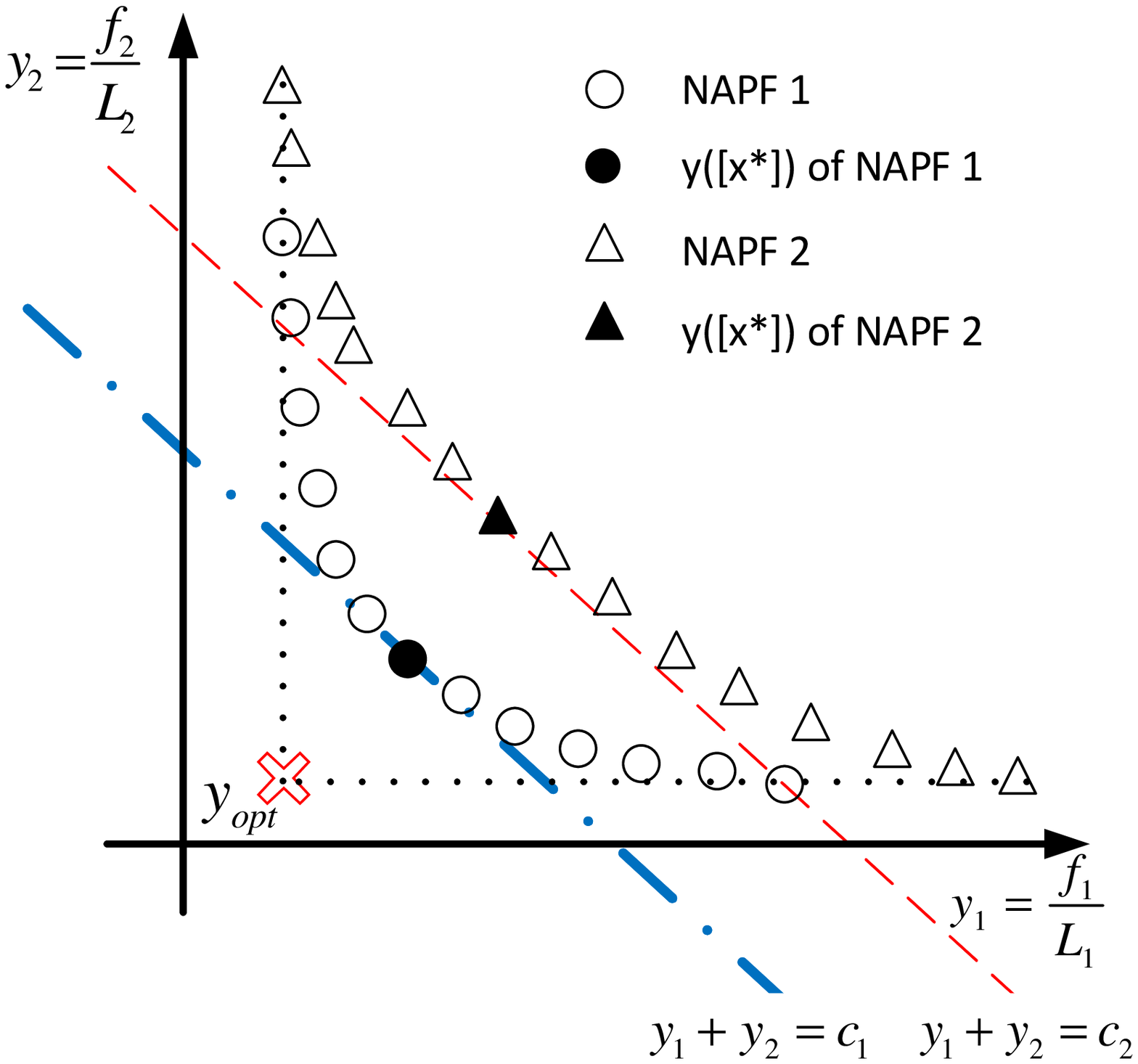} \\
    \mbox{(a)} & \mbox{(b)} \\
    \end{array}
\end{equation*}
  \caption{Geometric interpretation of the WS approach to MCDM. (a) The graph of a hyperplane in a 2D objective function space is a line  $y_1+y_2= c$ with intercept $c$. The $c_{min}^{WS}$
  is the minimum value of $c$ such that the line has nonempty intersection with the NAPF.
  (b) Two different shapes of NAPFs are considered. $c_1$ and $c_2$ represent the minimum values $c_{min}^{WS}$ associated with NAPF~1 and NAPF~2, respectively. A more bent NAPF can yield a smaller value of $c_{min}^{WS}$, i.e., $c_1<c_2$.}\label{fig_WS}
\end{figure*}

To begin with, we show that
the D\&C approach is equivalent to a WS approach that determines the final solution by adding weighting coefficients to objectives.

\begin{thm}\label{thm_knee_sol}
 The WS approach
\begin{equation}\label{eq_WS}
 \arg_{\bm{x}}\min_{ \bm{x}  \in \mathcal{PS}^{A} }    \sum_{n=1}^N   w_n f_n (\bm{x})
\end{equation}
  where
  \begin{equation}\label{eq_coe}
w_n=\frac{1}{L_n}
  \end{equation}
\end{thm}
 is equivalent to the  D\&C approach using the comparison rule in~(\ref{eq_pref2}).

{\emph{Proof:}
  This can be readily verified by noting that
\begin{equation*}
  \begin{split}
                   & [\bm{x}_j] \prec_{\bm{k}}  [\bm{x}_i] \\
\Leftrightarrow  { }  &  IP(\bm{x}_i \rightarrow \bm{x}_j)>0 \\
\Leftrightarrow  { }  &   \sum_{n=1}^N \frac{f_n(\bm{x}_i) - f_n(\bm{x}_j)    }{L_n} > 0 \\
\Leftrightarrow  { }  &  \sum_{n=1}^N \frac{f_n(\bm{x}_j) }{L_n} <  \sum_{n=1}^N \frac{f_n(\bm{x}_i) }{L_n}
   \end{split}
\end{equation*}
for all $\bm{x}_j, \bm{x}_i \in \mathcal{PS}^{A}$.
\hfill $\Box$

The WS approach assigns weighting coefficients
 to all objectives according to~(\ref{eq_coe}), and selects the solution that corresponds with the minimum sum.
It is different from conventional WS methods
in that it does not require preference inputs of the DM.
Define
\begin{equation*}
\begin{split}
  \bm{y} (\bm{x})  { }={ } &
   \left[
    \begin{array}{cccc}
     y_1(\bm{x}) & y_2(\bm{x}) & \ldots & y_N(\bm{x}) \\
    \end{array}
  \right]^T
   \\
  { }={ }   &
 \left[
    \begin{array}{cccc}
    \frac{ f_1(\bm{x})}{L_1} & \frac{ f_2(\bm{x})}{L_2} & \ldots & \frac{ f_N(\bm{x})}{L_N} \\
    \end{array}
  \right]^T.
\end{split}
\end{equation*}
 Because of the equivalence, the notation $[\bm{x}^*]$ adopted in the D\&C approach is also used here to denote the solution selected by the WS approach, i.e.,
\begin{equation*}
 \bm{x}^*= \arg_{\bm{x}}\min_{ \bm{x}  \in \mathcal{PS}^{A} }    \sum_{n=1}^N   \frac{f_n (\bm{x})}{L_n}.
\end{equation*}
From a geometric perspective,
the solution can be obtained by  moving a hyperplane
\begin{equation*}
y_1+y_2+...+y_N=c
\end{equation*}
 from a large $c$  in the direction  -$\bm{1}$  (-$\bm{1}$ represents the direction of decreasing the value of $c$)
 until the minimum value $c_{min}^{WS}$ ensuring nonempty intersection of the hyperplane and normalized APF (NAPF)
is achieved, as demonstrated in Fig.~\ref{fig_WS}(a). In this case, we have
\begin{equation*}
  c_{min}^{WS}= \min_{ \bm{x}  \in \mathcal{PS}^{A} }\sum_{n=1}^N   \frac{f_n (\bm{x})}{L_n}  .
\end{equation*}
For a smaller value of $c_{min}^{WS}$, the shape of the  NAPF can be more bent, as shown in Fig.~\ref{fig_WS}(b).

\begin{figure*}
  \centering
\begin{equation*}
    \begin{array}{cc}
     \includegraphics[width=9cm]{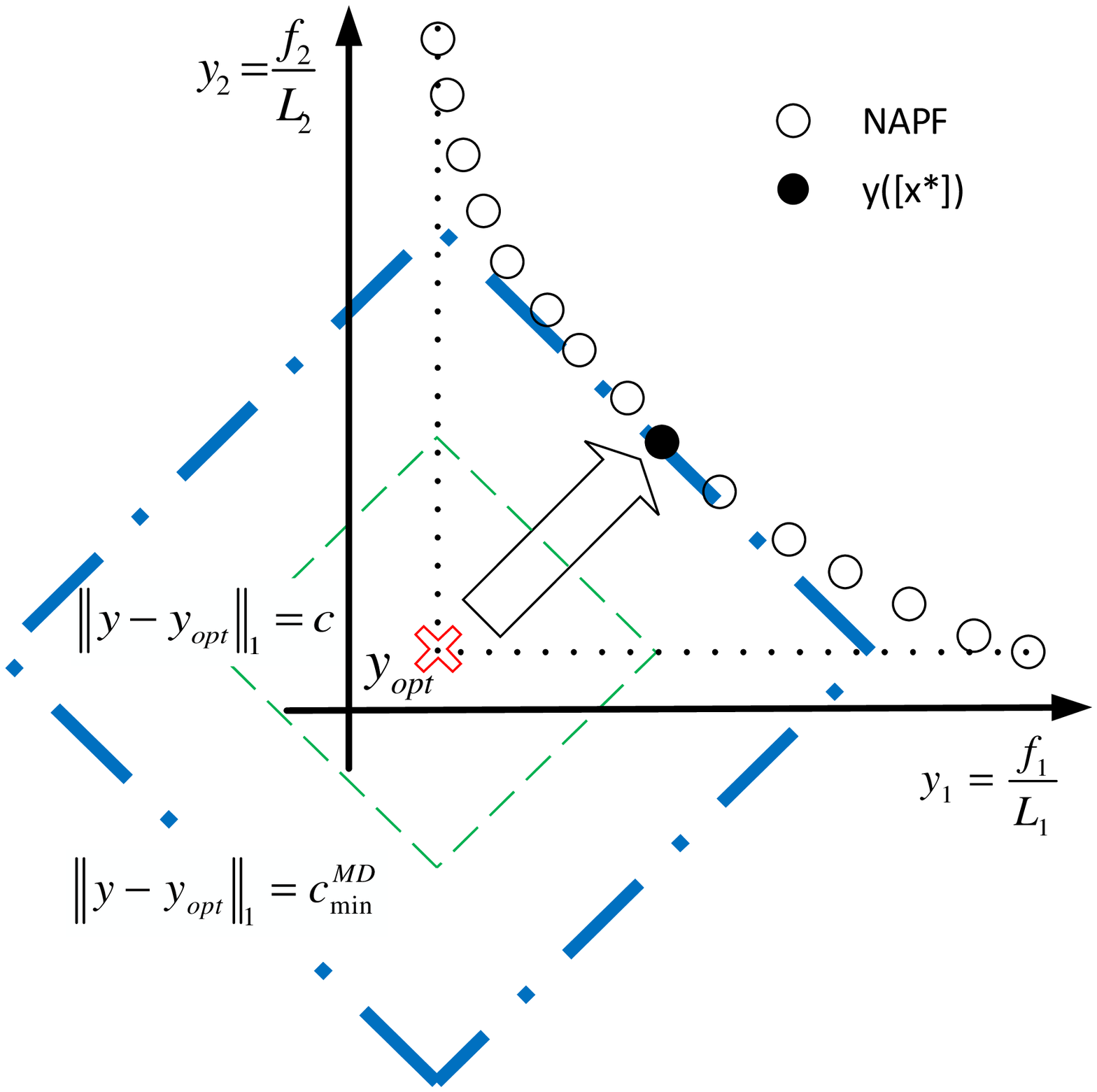} &  \includegraphics[width=9cm]{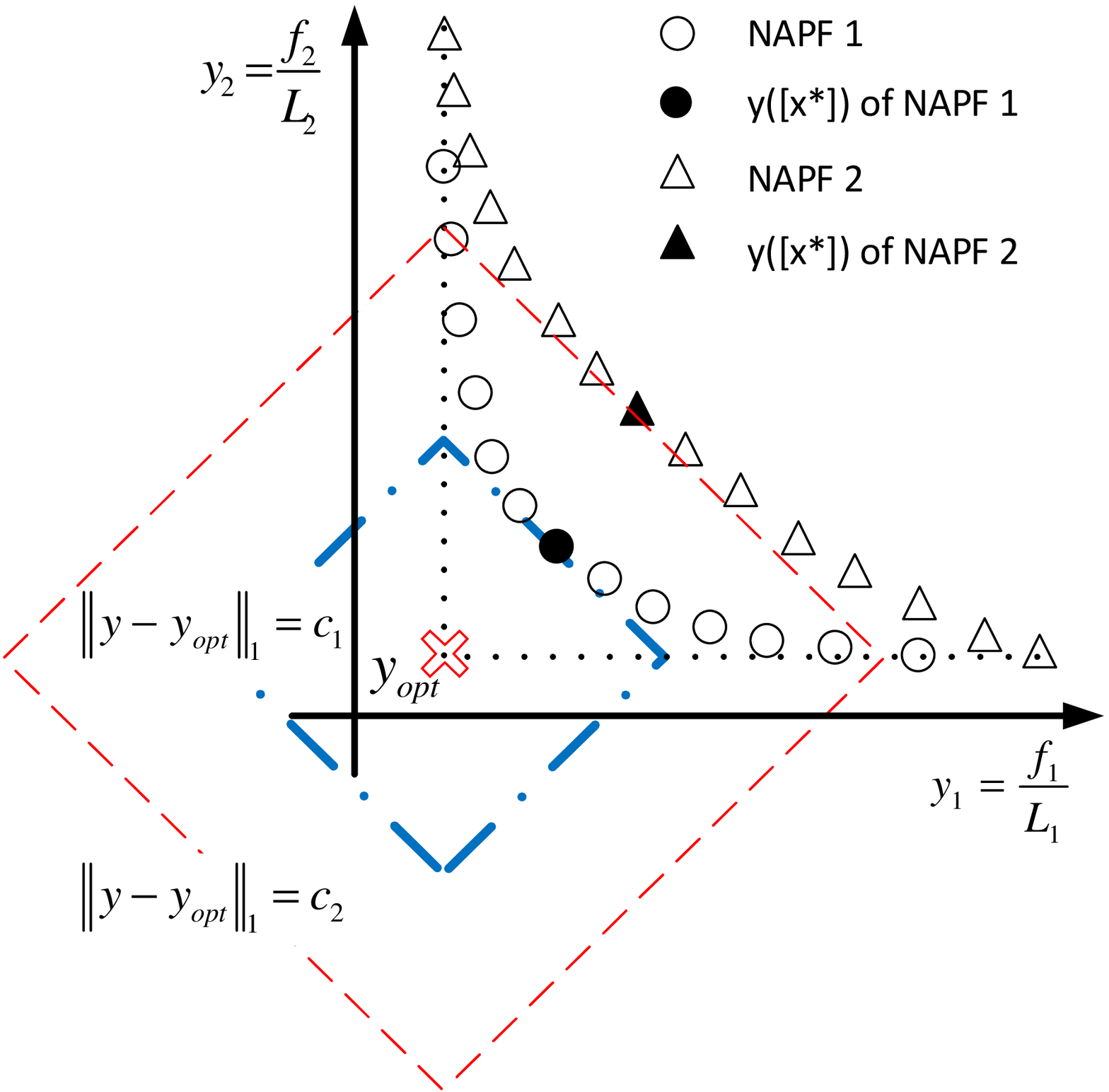} \\
    \mbox{(a)} & \mbox{(b)} \\
    \end{array}
\end{equation*}
  \caption{Geometric interpretation of the MMD approach to MCDM. (a) The graph
  of $||  \bm{y} - \bm{y}_{opt} ||_1= c$ in a 2D objective function space is a rhombus with center $\bm{y}_{opt}$ defined in~(\ref{eq_utopia}) and radius $c$, represented by the inner rhombus that has dashed edges. When $c=0$, the graph reduces to the point $\bm{y}_{opt}$.
  By enlarging $c$, the size of the rhombus increases, as indicated by the arrow.
   The minimum value of $c$  that yields nonempty intersection between the rhombus $||  \bm{y} - \bm{y}_{opt} ||_1= c$ and the NAPF is denoted by $c=c_{min}^{MMD}$, and
   the resulting rhombus is  $||  \bm{y} - \bm{y}_{opt} ||_1= c_{min}^{MMD}$, i.e., the outer rhombus that has dash-dot edges.
 This intersection is represented by $\bm{y}([\bm{x}^*])$ and  $[\bm{x}^*]$ is the solution set selected by the MMD approach.
  (b) Two different shapes of NAPFs are considered.   $c_1$ and $c_2$ represent the minimum values $c_{min}^{MMD}$ associated with NAPF~1 and NAPF~2, respectively.
  Since the shape of NAPF~1 is more bent than that of NAPF 2, we have $c_1<c_2$.}\label{fig_MMD}
\end{figure*}

With the help of Theorem~\ref{thm_knee_sol}, we can relate knee selection to the MMD.
Let
\begin{equation*}
  \ell_n=  \min_{ \bm{x}  \in \mathcal{PS}^{A} }   f_n(\bm{x})
\end{equation*}
be the minimum value in the $n$th dimension, and denote
  \begin{equation}\label{eq_utopia}
\bm{y}_{opt}  =
\left[
  \begin{array}{cccc}
    \frac{\ell_1}{L_1} & \frac{\ell_2}{L_2} & ... & \frac{\ell_N}{L_N} \\
  \end{array}
\right]^T
  \end{equation}
as the ideal vector after normalization.
The MMD approach to MCDM selects the point in the NAPF that is closest to the normalized ideal vector:
 \begin{equation}\label{eq_MMD}
   \bm{x}^*=   \arg_{\bm{x}}  \min_{ \bm{x}  \in \mathcal{PS}^{A} }  ||  \bm{y}(\bm{x}) - \bm{y}_{opt} ||_1
 \end{equation}
where $\bm{y}_{opt}$ is defined in~(\ref{eq_utopia}) and
   $||\cdot||_1$ represents the Manhattan norm (also termed 1-norm or taxicab norm), i.e., $||\bm{y}||_1=\sum_{n=1}^{N} | y_n|$.

The use of the Manhattan norm in the MMD approach described by~(\ref{eq_MMD}) establishes the connection with the weighted sum approach as shown in the following theorem.

\begin{thm}\label{thm_knee_sol3}
 The MMD approach is equivalent to the WS approach. In other words, we have
  \begin{equation}\label{eq_WS_equi}
  \min_{ \bm{x}  \in \mathcal{PS}^{A} }  ||  \bm{y}(\bm{x}) - \bm{y}_{opt} ||_1  \sim \min_{ \bm{x}  \in \mathcal{PS}^{A} }    \sum_{n=1}^N   \frac{f_n (\bm{x})}{L_n}
  \end{equation}
where ``$\sim$'' denotes ``equivalent to.''
\end{thm}

{\emph{Proof:}
 This can be readily verified by noting that
 \begin{equation*}
 \begin{split}
    &
   ||  \bm{y}(\bm{x}) - \bm{y}_{opt} ||_1 = \sum_{n=1}^N |\frac{f_n(\bm{x}) -\ell_n }{L_n}| =   \sum_{n=1}^N  \{ \frac{f_n(\bm{x})  }{L_n}   -  \frac{\ell_n }{L_n} \} \\
   ={ }  &    \sum_{n=1}^N   \frac{f_n(\bm{x})  }{L_n}   -  \sum_{n=1}^N  \frac{\ell_n }{L_n}
 \end{split}
 \end{equation*}
 where the second equality comes from the fact that  $f_n(\bm{x})  \geq \ell_n$ for all $\bm{x}  \in \mathcal{PS}^{A}$. Therefore, we have
   \begin{equation*}
    \begin{split}
 &  \min_{ \bm{x}  \in \mathcal{PS}^{A} }  ||  \bm{y}(\bm{x}) - \bm{y}_{opt} ||_1  \sim \min_{ \bm{x}  \in \mathcal{PS}^{A} }  \{   \sum_{n=1}^N  \frac{f_n (\bm{x})}{L_n}    -  \sum_{n=1}^N  \frac{\ell_n }{L_n}   \}  \\
   \sim { } & \min_{ \bm{x}  \in \mathcal{PS}^{A} }\sum_{n=1}^N   \frac{f_n (\bm{x})}{L_n}
     \end{split}
  \end{equation*}
  because the term $ \sum_{n=1}^N  \ell_n /L_n$ is a constant.
\hfill $\Box$

\begin{figure*}
  \centering
\begin{equation*}
\begin{array}{ccc}
\hspace{-0.5cm}  \includegraphics[width=6cm]{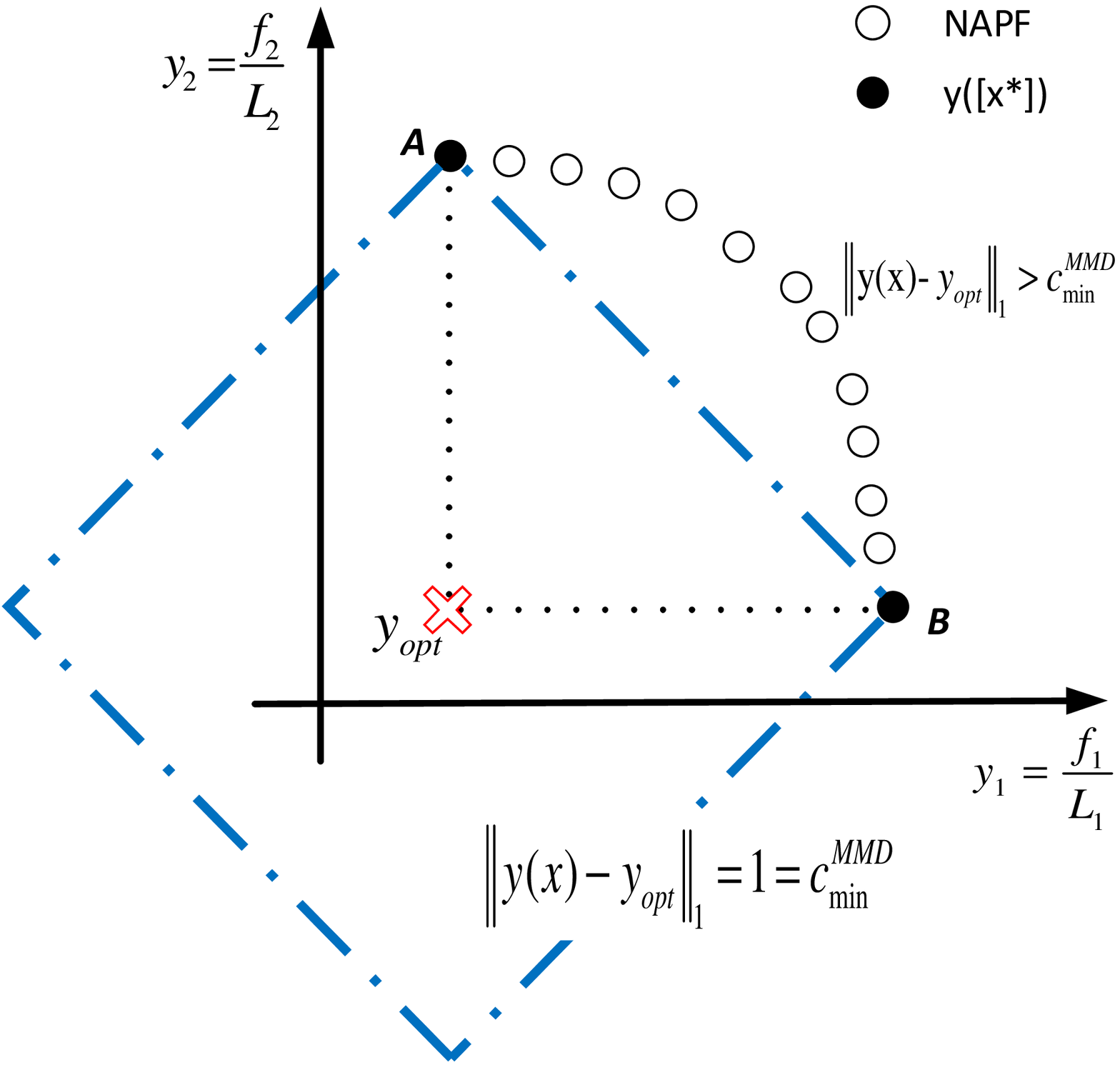} & \includegraphics[width=6cm]{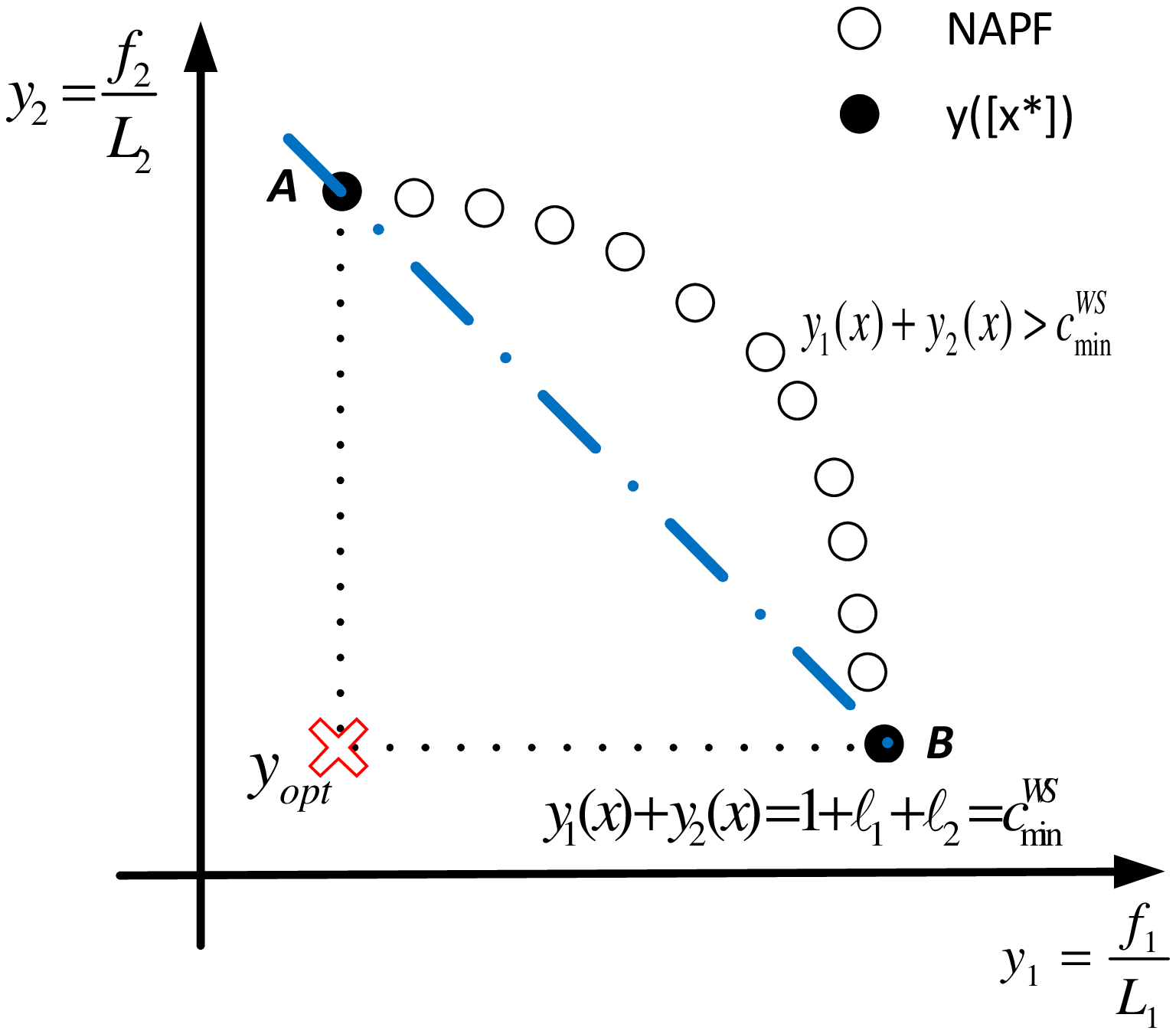} &  \includegraphics[width=6cm]{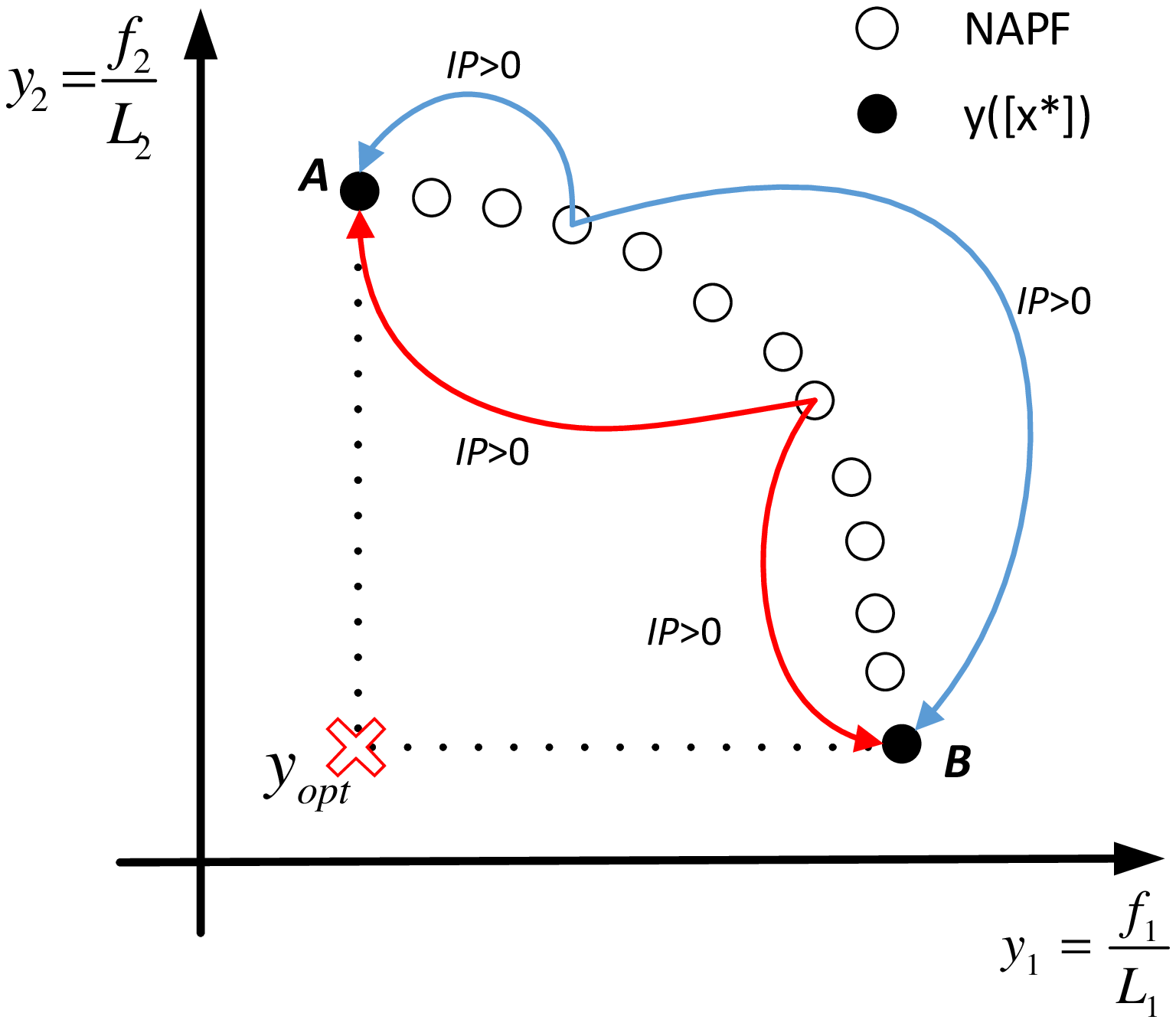}\\
  \mbox{(a)} & \mbox{(b)} & \mbox{(c)} \\
\end{array}
\end{equation*}
  \caption{The knee $\bm{y}([\bm{x}^*])$ contains two points, i.e., $\bm{y}([\bm{x}^*])=\{  \bm{A},\bm{B} \}$. (a) Nondominated vectors outside the rhombus yield $||  \bm{y}(\bm{x}) - \bm{y}_{opt} ||_1>c_{min}^{MMD}$.
(b) Nondominated vectors located on the right-hand side of the line $ y_1(\bm{x})+ y_2(\bm{x}) = c_{min}^{WS} $ yield $ y_1(\bm{x})+ y_2(\bm{x}) > c_{min}^{WS} $.
(c) Transition from middle nondominated vectors to end vectors is preferred because a positive net  improvement percentage, denoted by $IP>0$, can be achieved.
 }\label{fig_equi_class}
\end{figure*}

From an algebraic perspective, two observations can be made  from Theorem~\ref{thm_knee_sol3}.
First, the MMD approach is efficient because
evaluating the Manhattan norm can be realized by efficient matrix calculations:
\begin{equation*}
    ||  \bm{y}(\bm{x}) - \bm{y}_{opt} ||_1 = \sum_{n=1}^N  \{ \frac{f_n(\bm{x})  }{L_n}   -  \frac{\ell_n }{L_n} \} =\bm{1}^T (  \bm{y}(\bm{x}) - \bm{y}_{opt})
\end{equation*}
where
$\bm{1}$ represents the vector with all-one entries.
Second, although the MMD and WS approaches are equivalent,
the MMD approach is generally preferred.
When a term $f_n(\bm{x}) / L_n$ in~(\ref{eq_WS}) is too large compared to the remaining terms,
  the WS approach neglects the remaining ones, which may yield difficulty searching for the final solution.
  Note that this difficulty cannot be avoided by simply normalizing objectives.
This is because any normalizing constant $\alpha_n$ in the $n$th dimension will enter the maximum spread so that the normalizing effect is cancelled in the ratio, i.e.,
\begin{equation*}
\frac{\tilde{f}_n}{\tilde{L}_n}   =\frac{(f_n/\alpha_n)}{(L_n/\alpha_n)}=\frac{f_n}{L_n}
\end{equation*}
where $\tilde{f}_n=f_n/\alpha_n$ and $\tilde{L}_n=L_n/\alpha_n$ represent the normalized objective and associated maximum spread, respectively.
By contrast,
all the terms $(f_n(\bm{x}) -\ell_n) / L_n$ belong to $[0,1]$
in the MMD approach, avoiding the problem of one term dominating the remaining terms.

For a geometric interpretation, we consider
the graph of
\begin{equation*}
   ||  \bm{y} - \bm{y}_{opt} ||_1= c
\end{equation*}
 which in $\mathbb{R}^2$ is a rhombus with center $\bm{y}_{opt}$  and radius $c$.
The selected  $\bm{y}([\bm{x}^*])$ can be obtained by gradually enlarging $c$ until the graph intersects normalized $\mathcal{PF}^{A}$.
The value
\begin{equation*}
  c_{min}^{MMD}=  \min_{ \bm{x}  \in \mathcal{PS}^{A} }  ||  \bm{y}(\bm{x}) - \bm{y}_{opt} ||_1
\end{equation*}
 is the minimal value for nonempty intersection of the rhombus and normalized $\mathcal{PF}^{A}$, as explained in Fig.~\ref{fig_MMD}(a).
  In addition to indicating how close the selected $\bm{y}([\bm{x}^*])$ is to the ideal vector~$\bm{y}_{opt}$,
the value of~$c_{min}^{MMD}$ may reveal how objectives affect each other.
Similarly to the role of~$c_{min}^{WS}$, a smaller value of $c_{min}^{MMD}$ implies a more bent NAPF, as illustrated in Fig.~\ref{fig_MMD}(b).

Theorems~\ref{thm_knee_sol} and~\ref{thm_knee_sol3} provide
the proposed MMD approach with rich geometric and algebraic interpretations because of the overall equivalence established.
In summary, from a geometric perspective,
 the proposed approach selects $\bm{x}^*$ that has the MMD from the normalized ideal vector~$\bm{y}_{opt}$;
$\bm{x}^*$ can also be obtained by either enlarging the radius of the rhombus with the normalized ideal vector as the center
or moving the hyperplane toward the direction of decreasing the value of intercept until nonempty intersection with the NAPF cannot be achieved;
and the approach is equivalent to the knee selection method described by the D\&C approach.
Computationally, the MMD approach is more efficient and elegant than
the D\&C approach that requires iterations of pairwise comparisons to yield $\bm{x}^*$;
it is more effective than the WS approach in certain situations in which the WS approach has difficulty searching for the final solution; and the approach can be considered as a systematic way to assign weighting coefficients to objectives,
which is generally difficult when a large number of objectives are involved.

Visualizing knee selection has been examined solely in the case where the APF is convex.
We examine the knee selection for other shapes of fronts.
Fig.~\ref{fig_equi_class} shows a concave APF in a 2-D objective function space.   Suppose that
\begin{equation*}
{\small
  \bm{A}=
\left[
  \begin{array}{cc}
    A_1 & A_2 \\
  \end{array}
\right]^T=  \bm{y}(\bm{x}_1)
\mbox{ and }
\bm{B}=
\left[
  \begin{array}{cc}
    B_1 & B_2 \\
  \end{array}
\right]^T=  \bm{y}(\bm{x}_2)
}
\end{equation*}
 are the two extreme vectors.
According to~(\ref{eq_utopia}), we have
\begin{equation*}
 \bm{y}_{opt}=
 \left[
  \begin{array}{cc}
    A_1 & B_2 \\
  \end{array}
\right]^T=
 \left[
  \begin{array}{cc}
    \frac{\ell_1}{L_1} & \frac{\ell_2}{L_2} \\
  \end{array}
\right]^T.
\end{equation*}
 For the MMD approach, we have
\begin{equation}\label{eq_concave_MMD}
\begin{split}
            &  ||  \bm{y}(\bm{x}) - \bm{y}_{opt} ||_1= \frac{f_1(\bm{x}) - \ell_1    }{L_1}  +\frac{f_2(\bm{x}) - \ell_2    }{L_2}   \\
 { } ={ }   &
 \left\{
   \begin{array}{ll}
     0+1=1, & \hbox{ if } \bm{x}=\bm{x}_1 \\
     1+0=1, & \hbox{ if } \bm{x}=\bm{x}_2
   \end{array}
 \right..
 \end{split}
\end{equation}
 Other vectors in the APF yield $||  \bm{y}(\bm{x}) - \bm{y}_{opt} ||_1>1$ and thus
$\bm{y}([\bm{x}^*])=\{  \bm{A}, \bm{B} \}$.
By rearranging terms in~(\ref{eq_concave_MMD}), we have
\begin{equation*}
 y_1(\bm{x})+ y_2(\bm{x}) =    \frac{f_1(\bm{x})}{L_1} + \frac{f_2(\bm{x})  }{L_2}
 =  1 + \ell_1 + \ell_2
\end{equation*}
for  $\bm{x}\in [\bm{x}^*]$.
Therefore, when the WS approach is used,  points $\bm{A}$ and $\bm{B}$ are on the same line, leading to the minimum weighted sum.
We see that  the D\&C approach produces the same result:
transition from extreme point $\bm{A}$ or $\bm{B}$ to middle points is not allowed because the
 improvement percentage in one dimension is less than the  degradation percentage in the other dimension.

Selection of extreme vectors when the shape of an APF is a concave curve has further implication.
In a 2-D space, if the shape of an APF can be represented by a line segment,
then all vectors on the line will be selected by the proposed methodology.
This is because  a line can be regarded as a degenerate case of a concave curve so that
 all vectors on the line become extreme vectors.
Since a line in a 2-D space generalizes to a plane in a 3-D space,
if the shape of an APF can be contained within a plane in a 3-D space,
 then vectors on the plane will be chosen. After normalization, the chosen vectors yield the same MMD from the ideal vector.
 For any APFs represented  by concave surfaces in a 3-D space,
extreme vectors on the concave surfaces are to be selected.
 This can be understood by noting that a concave curve in a 2-D space generalizes to a concave surface in a 3-D space.

Finally, we discuss how a vector is selected by the MMD approach when a discontinuous front in a 2-D space is involved.
After normalizing all vectors in the front, we can construct a line segment by connecting the two extreme vectors, i.e.,
the vector with the smallest value in $y_1$ and the vector with the smallest value in $y_2$.
Similar to the scenario considered in Fig.~\ref{fig_equi_class}(a),
the following observations can be made: normalized vectors $\bm{y}$ on the line segment,
on the right-hand side of the line segment, and on the left-hand side of the line segment
yield $||  \bm{y}- \bm{y}_{opt} ||_1=1$, $||  \bm{y}- \bm{y}_{opt} ||_1>1$, and  $||  \bm{y}- \bm{y}_{opt} ||_1<1$, respectively.
 Therefore, the MMD approach selects a vector that is on the left-hand side of the line segment. If all the non-extreme vectors
 are on the right-hand side of the line segment, then extreme vectors are to be selected.

\begin{rmk}
Compromise programming to MCDM determines the final solution that is associated with the least distance from an ideal vector~\cite{pohekar2004CP,Zeleny82CP}.
The distance is measured in terms of $p$-norms combined with weighting coefficients prescribed by the DM.
The distance function in compromise programming reduces to the Manhattan distance used in our approach
when $p=1$ is assigned and  equal  weighting coefficients  are adopted.
Therefore, to some extent our analysis has connects not only knee selection  with WS methods, but also
WS methods with compromise programming.
\end{rmk}

\begin{figure}
  \centering
\begin{equation*}
\begin{array}{c}
\includegraphics[width=7cm]{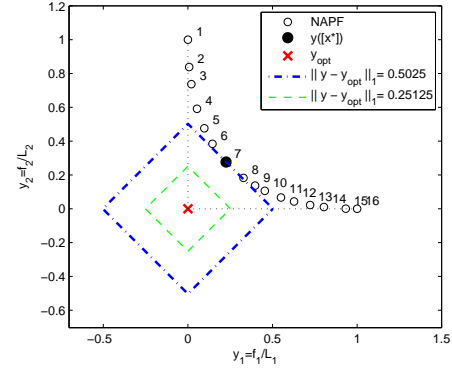} \\
  \mbox{(a)} \\
  \includegraphics[width=7cm]{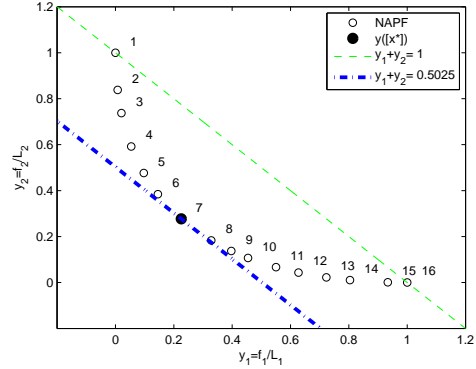} \\
\mbox{(b)}\\
\includegraphics[width=6cm]{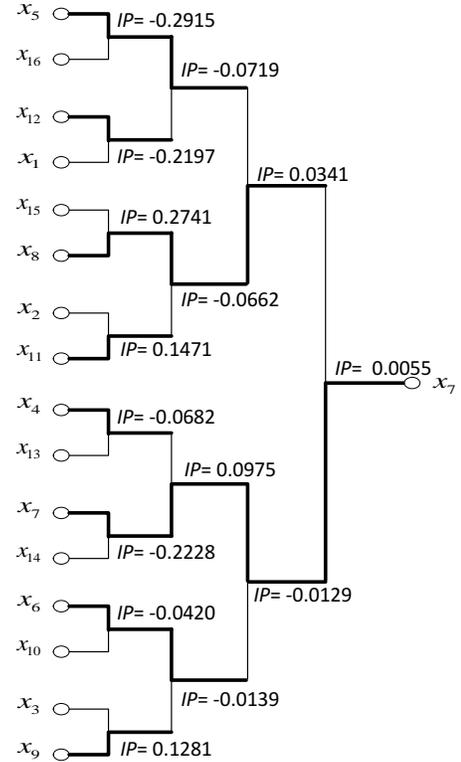}\\
\mbox{(c)} \\
\end{array}
\end{equation*}
  \caption{MCDM in MOP1 by (a) MMD approach; (b) WS approach; and (c) D\&C approach. Solutions are labeled based on the associated objective values of $f_1$.
In (c), $IP=IP(\bm{x}_i  \rightarrow \bm{x}_j )$ where $\bm{x}_i$ and $\bm{x}_j$ represent the upper and lower solutions in a pairwise comparison, respectively.
 }\label{fig_sim_2D}
\end{figure}

\begin{figure}
  \centering
\begin{equation*}
\begin{array}{c}
\includegraphics[width=7cm]{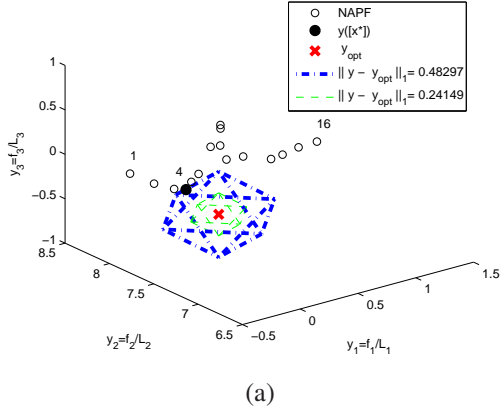} \\
  \mbox{(a)} \\
  \includegraphics[width=7cm]{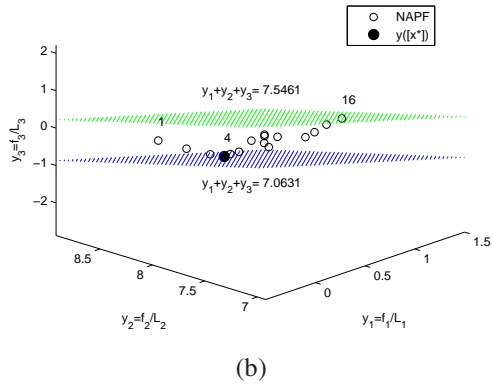} \\
\mbox{(b)}\\
\includegraphics[width=6cm]{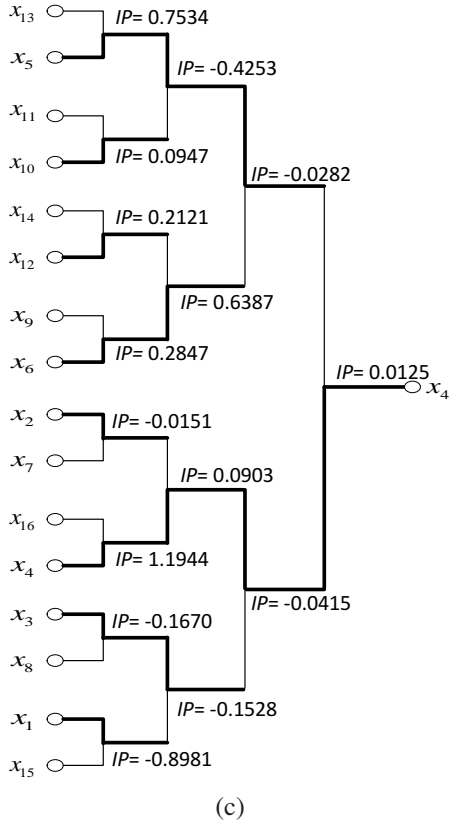}\\
\mbox{(c)} \\
\end{array}
\end{equation*}
  \caption{MCDM in MOP5 by (a) MMD approach; (b) WS approach; and (c) D\&C approach.  Solutions are labeled based on the associated objective values of $f_1$. In (c), $IP=IP(\bm{x}_i  \rightarrow \bm{x}_j )$ where $\bm{x}_i$ and $\bm{x}_j$ represent the upper and lower solutions in a pairwise comparison, respectively.}\label{fig_sim_3D}
\end{figure}

\section{Numerical Results}\label{sec_sim}

In this section, we examine various MCDM problems derived from MOPs to illustrate the proposed methodology.
The section is divided into three subsections.
The established equivalence is examined using 2-D, 3-D, and 5-D APFs in Section~\ref{subsec_equi}.
In Section~\ref{subsec_R_E}, practical concerns about these equivalent approaches are investigated.
In Section~\ref{subsec_explore}, several benchmark MOPs are employed to demonstrate the effectiveness of the MMD approach.
Finally, an MCDM problem with real-world data is considered in Section~\ref{subsec_data}.
All simulations have been performed using a desktop with Intel i7-4770, 3.40 GHz CPU, and 3.16 GB RAM.

\begin{figure*}
  \centering
\begin{equation*}
\begin{array}{cccc}
\includegraphics[width=4cm]{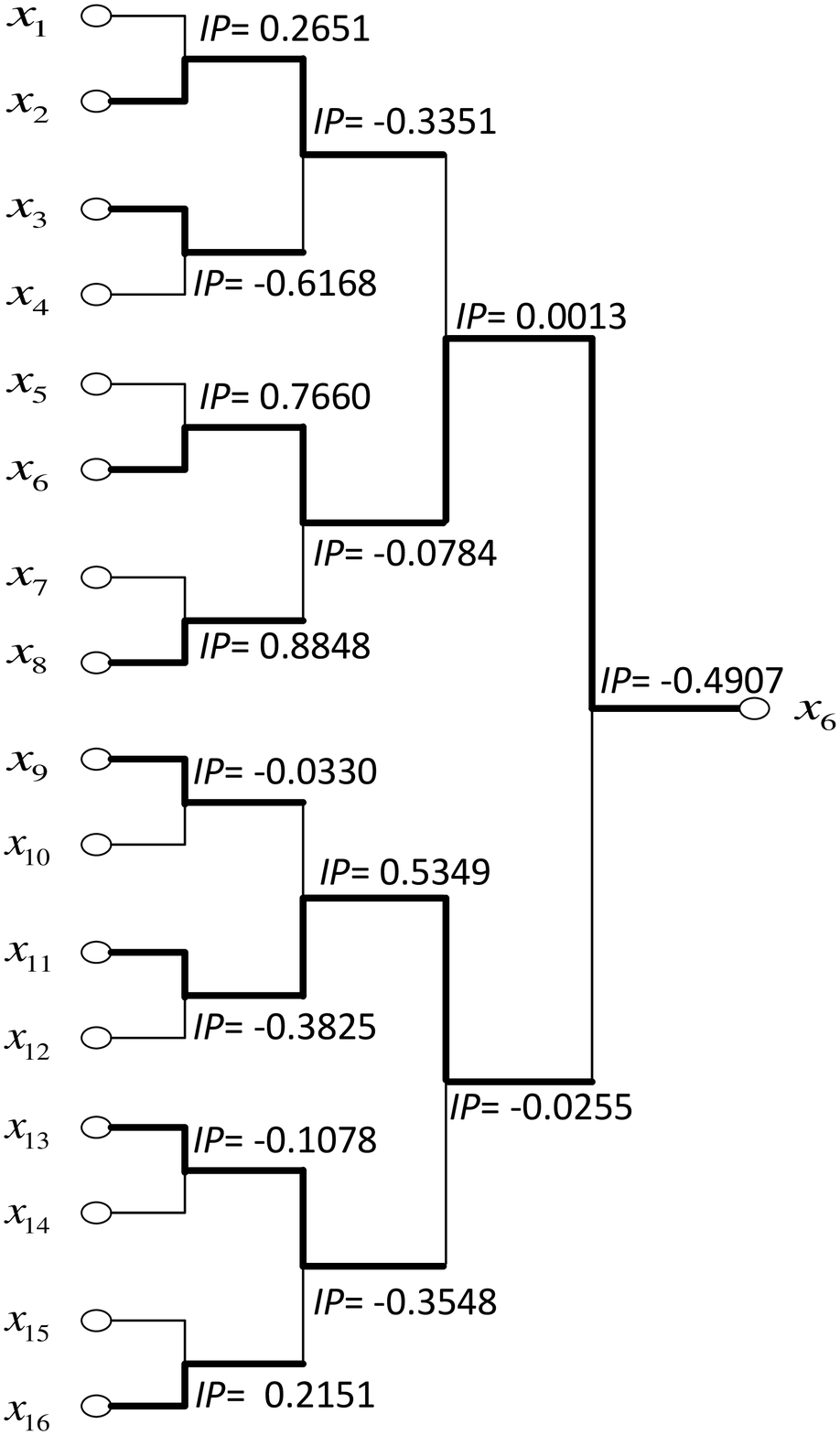} &  \includegraphics[width=4cm]{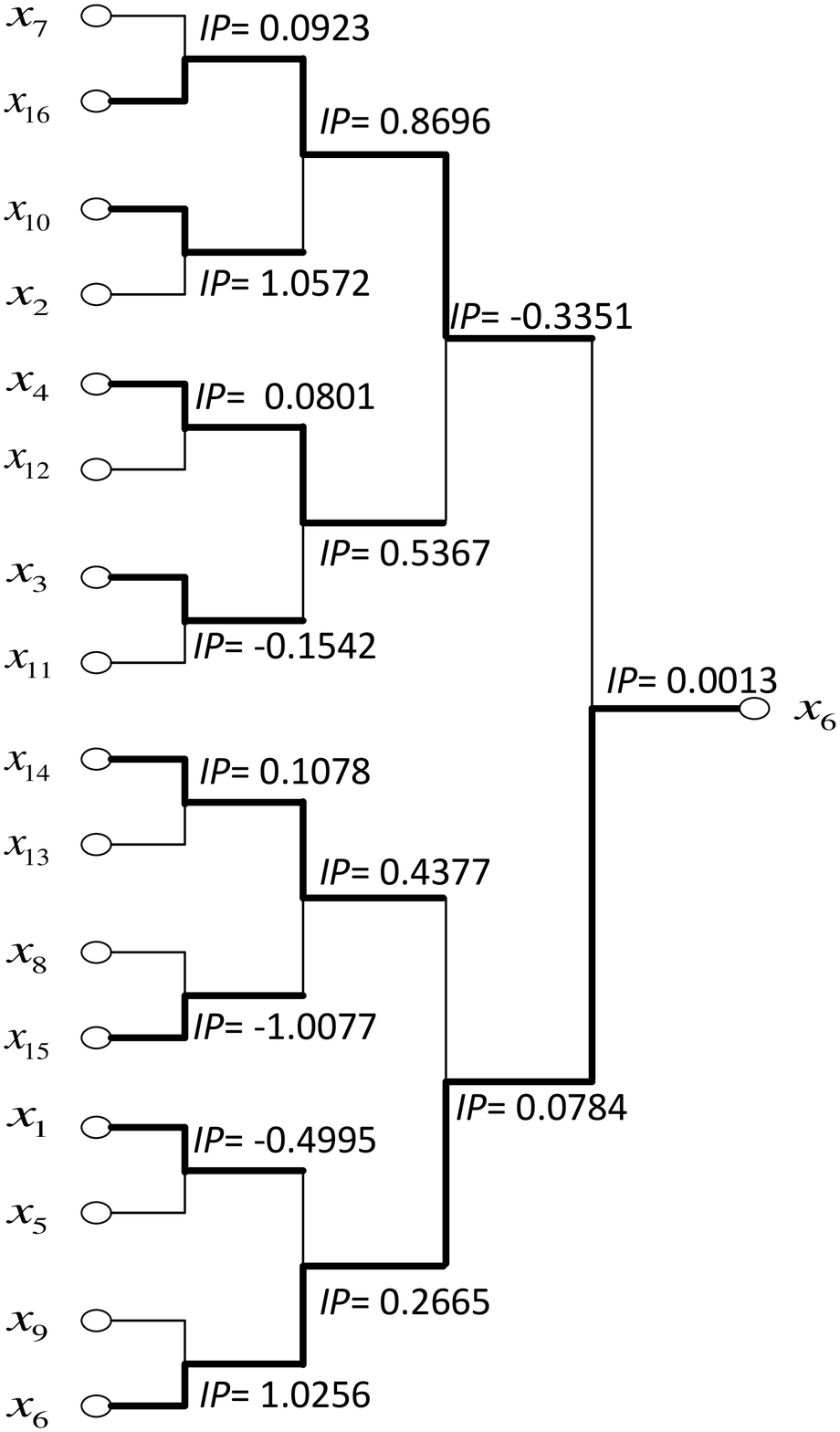}  & \includegraphics[width=4cm]{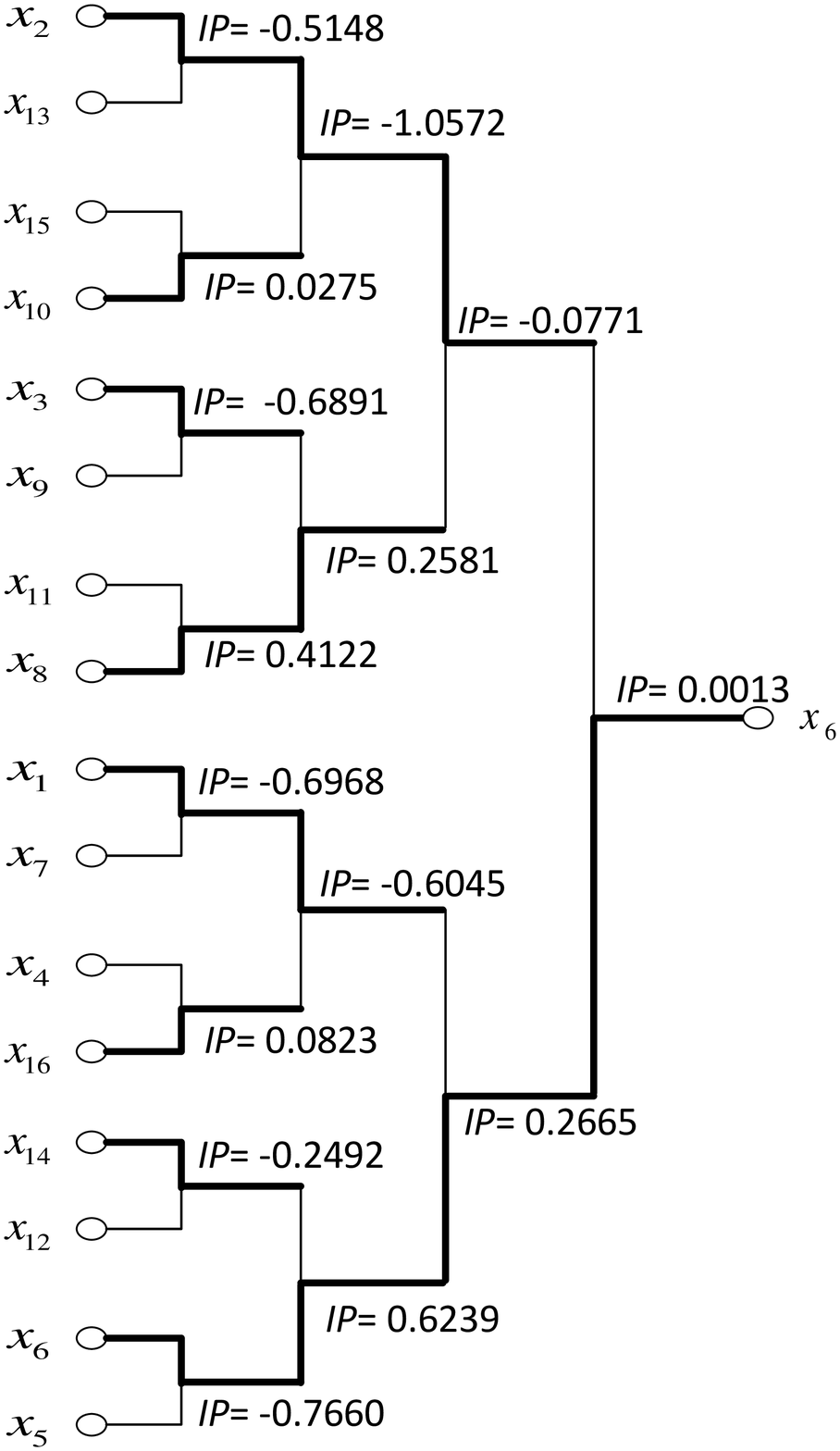} & \includegraphics[width=4cm]{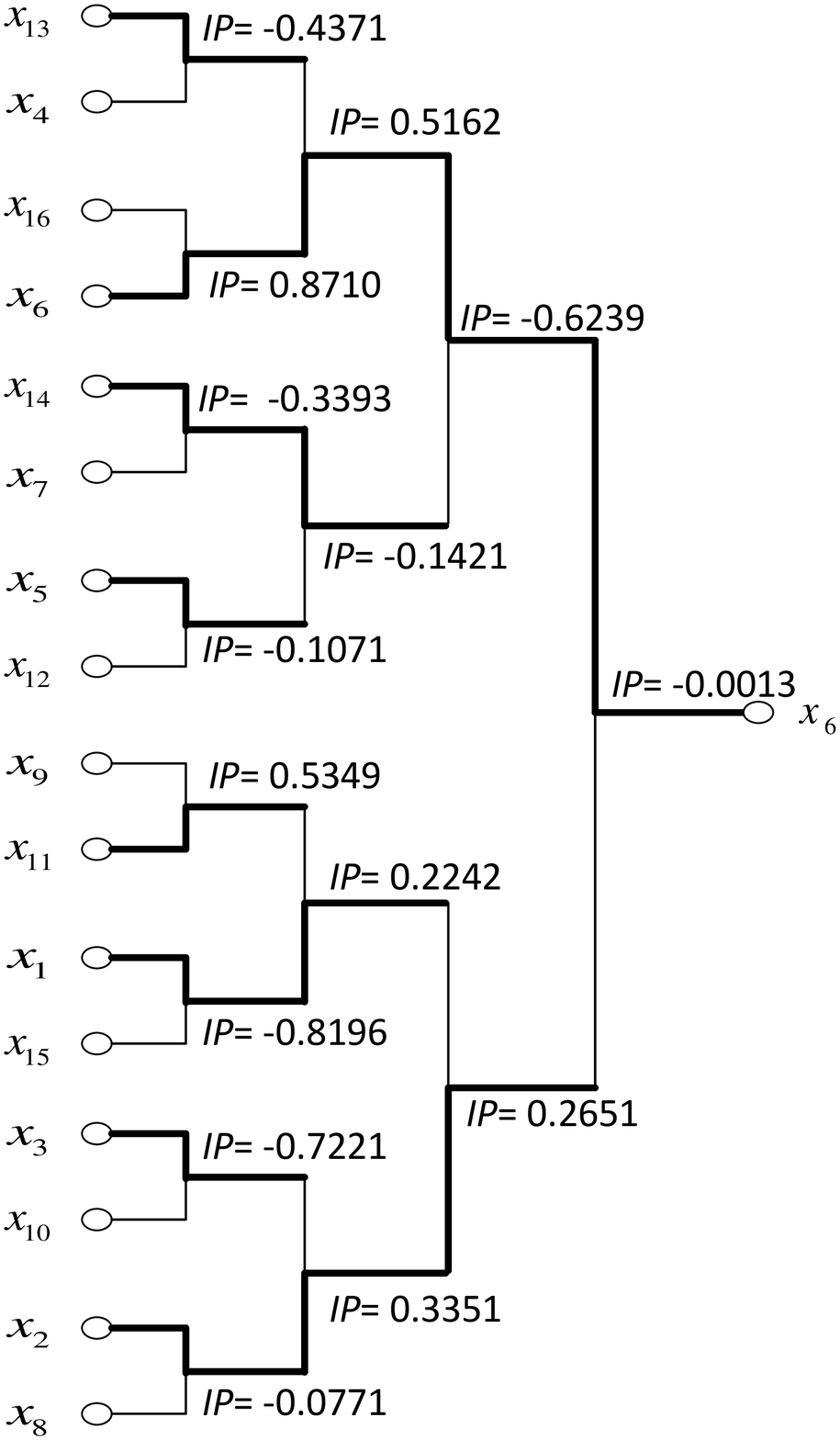}\\
  \mbox{(a)} & \mbox{(b)} & \mbox{(c)} & \mbox{(d)} \\
\end{array}
\end{equation*}
  \caption{Four random trials of pairwise comparisons using the D\&C approach to MCDM in DTLZ1. Different comparing orders in (a)--(d) lead to the same solution, where
 $IP=IP(\bm{x}_i  \rightarrow \bm{x}_j )$ with $\bm{x}_i$ and $\bm{x}_j$  representing the upper and lower solutions in a pairwise comparison, respectively.}\label{fig_sim_5D}
\end{figure*}

\subsection{Equivalence Analysis}\label{subsec_equi}

The MOP1, MOP5, and DTLZ1~\cite{MOEA_bk1} were solved by
the multiobjective artificial immune algorithm in~\cite{TSG_15} to obtain APFs.
For an illustrative purpose, a small population size of 16 was adopted.
For a 2-D illustration, i.e., $N=2$,
we considered the MCDM in MOP1.
As shown in Figs.~\ref{fig_sim_2D}(a) and~\ref{fig_sim_2D}(b),
the MMD and WS approaches are associated with a rhombus and a line, respectively.
 In~Fig.~\ref{fig_sim_2D}(a),   the rhombus $||  \bm{y} - \bm{y}_{opt} ||_1 =c$ reduces to the point $\bm{y}_{opt}$ when $c=0$.
  The inner rhombus with dashed edges has a value of $c=0.33443$, yielding empty intersection with the NAPF. By enlarging the value of $c$, the rhombus size increases.
  The outer rhombus with dash-dot edges that has a value of $c=0.66885$ is obtained by enlarging $c$ from $c=0$ until nonempty intersection with the NAPF is achieved.
  This intersection is represented by $\bm{y}(\bm{x}_7)$ and, therefore, the MMD selects the point $\bm{x}_7$.
  It should be noted that although we used a geometric interpretation to realize the MMD approach, the algebraic formula in~(\ref{eq_MMD}) should be used in practice.
Fig.~\ref{fig_sim_2D}(c) shows a random order of pairwise comparisons based on the D\&C approach.
All three approaches yield the same result, as proven in our analysis of their equivalence.
The vector associated with the final solution~$\bm{x}_7$ geometrically lies in the knee region of the APF, corresponding with knee selection
 and coinciding with our geometric intuition for a knee.

For a 3-D case, i.e., $N=3$, a line and a rhombus for the WS and MMD approaches become a plane and a regular octahedron, respectively.
Fig.~\ref{fig_sim_3D} presents the MCDM in MOP5.
In Fig.~\ref{fig_sim_3D}(a),
the MMD approach is interpreted as enlarging the radius of a regular octahedron centered at the ideal vector $\bm{y}_{opt}$ until nonempty intersection with the NAPF is achieved.
In Fig.~\ref{fig_sim_3D}(b),
the WS approach is interpreted as moving a plane in the direction of its normal vector -$\bm{1}$ while intersection with the NAPF must be ensured, leading to the minimum value of the weighted sum.
In Fig.~\ref{fig_sim_3D}(c), the D\&C approach with a random order of pairwise comparisons is applied, producing the same solution as the MMD and WS approaches do.

For a 5-D scenario,
the MCDM in scalable DTLZ1 is considered.
While geometric visualization becomes impossible, our algebraic formulas for MCDM can still be applied.
The MMD and WS approaches yield the same solution, bold marked in Table~\ref{tab_sim_5D}.
Four random trials of pairwise comparisons using the D\&C approach
 are performed to demonstrate that the approach is independent of the comparing order, shown in Fig~\ref{fig_sim_5D}.

\subsection{Practical Concerns}\label{subsec_R_E}

We showed that the MMD, WS, and D\&C approaches yielded the same final solutions in 2-D, 3-D, and 5-D MCDM problems.
Although these approaches are theoretically equivalent, in practice there are some situations in which the WS approach can have difficulty searching for the final solution
and the D\&C approach can consume relatively more computational time.

In Fig.~\ref{fig_WS_fail}(a), the values in $f_1$ are much larger than the associated maximum spread $L_1$,
but the differences between the values in $f_2$ and the associated maximum spread $L_2$  are relatively small.
We have $f_1/L_1\gg f_2/L_2$ and hence, the term $f_1/L_1$ dominates the term $f_2/L_2$ in the weighted sum.
As shown in Table~\ref{tab_WS_fail} and Fig.~\ref{fig_WS_fail}(b), all solutions have almost the same weighted sum, but the MMD and D\&C approaches can readily distinguish among the solutions. In this situation the WS approach has difficulty in finding the final solution.

To evaluate the corresponding computational time,
we examine the MCDM in DTLZ1, DTLZ2, MOP1--7, MOP-C1 Binh, MOP-C1 Osyczka, MOP-C1 Viennet, MOP-C1 Tanaka, and ZDT1--3~\cite{MOEA_bk1}.
Larger population sizes are adopted for statistical analysis.
 Table~\ref{tab_sim_time} summarizes the comparisons.
To facilitate ensuing discussions, we label three groups of simulation results as category~1 (C1),
category~2 (C2), and category~3 (C3). In C1 comparisons, 3000 simulation runs seem to allow for
relatively stable evaluation of average computational time. Among these comparisons,
the D\&C approach consume more computational time than the other approaches.
The population size of the employed MOEA is related to the number of solutions (or problem size) in the MCDM process,
and a larger size implies more computational efforts.
C2 comparisons illustrate that  computational time of the D\&C approach
increases more rapidly than the MMD and WS approaches upon increasing the problem size.
This is because the D\&C approach must perform pairwise comparisons iteratively and the number of comparisons is directly related to the problem size.
Various MCDM problems are examined in C3 comparisons, demonstrating that
the MMD and WS approaches are more computationally efficient than the
D\&C approach.

\subsection{Further Exploration}\label{subsec_explore}

For the purpose of a better understanding,
the MMD approach is applied to the MCDM in
commonly used ZDT and DTLZ test suites.
Since a few of these MOPs yield the same PFs, they are combined.
In addition, MOP4 and MOP6 from~\cite{MOEA_bk1} are included for comparison.
Fig.~\ref{fig_ZDT_DTLZ} presents the results in which the MCDM is performed on population sampled from the true PFs.
Normalized samples closest to the ideal vector in the sense of the Manhattan distance are to be selected.
For convex shapes of PFs in ZDT1 and ZDT4,
samples located in the knee region are selected as expected.
Because the PF in ZDT2 and ZDT6 (2-D problems) has the shape of a concave curve
and that in DTLZ2--4 (3-D problems) has the shape of a concave surface,
 extreme samples in each dimension are selected.
Since samples of the PF in DTLZ1 are contained within a plane,
 all of them are chosen and considered as equivalent.

For discontinuous PFs, we refer to ZDT3, DTLZ7, MOP4, and MOP6.
It is informative to compare ZDT3, MOP4, and MOP6.
Consider the line segment that connects the extreme vectors.
 For ZDT3, most samples are on the left-hand side of the line, and the one that is most distant from the line in the sense of the Manhattan norm is selected.
 By contrast,  MOP4 and MOP6  have most samples on the right-hand side of the line; however, the selected samples are an exception that is on the left-hand side but close to the line.

\subsection{Real-World Application}\label{subsec_data}

The MMD approach is further used to solve a real-world MCDM problem about a future plant layout of a leading IC packaging company in Taiwan~\cite{kuo2008use}.
It is desired that the plant layout can have certain features
measured by the flow distance ($f_1$), adjacency score ($-f_2$ where the negative sign indicates a larger-the-better quantity), shape ratio ($f_3$), flexibility ($-f_4$), accessibility ($-f_5$), and maintenance ($-f_6$). In this problem, there are 18 layout alternatives ($\bm{x}_1$--$\bm{x}_{18}$), generated by a commercial software program termed Spiral.
  Existing MCDM approaches are included for comparison:
  grey relational analysis  (GRA),   data envelopment analysis (DEA),\footnote{The DEA approach is further combined with an analytical hierarchy process.}
  the technique for order preference by similarity to an ideal solution (TOPSIS), and simple additive weighting (SAW) ~\cite{kuo2008use,yang03hierarchical,yang2007multiple}.

In practice a DM uses various analysis tools, compares the results, and then selects the final alternative when addressing an MCDM problem.
Table~\ref{tab_RealData} presents the ranking of the alternatives.\footnote{Due to space consideration, only the top 10 alternatives are listed.}
It is worth mentioning that while state-of-the-art MCDM approaches have different mechanisms,
most of them put alternatives $\bm{x}_{11},\bm{x}_{15}$, and $\bm{x}_{17}$ in the top-3 list.
Alternative $\bm{x}_{15}$ should be selected because it has the top ranking among most MCDM approaches.
The proposed MMD approach is consistent with this selection.

Although most existing MCDM approaches yield the same final result,
the proposed MMD approach is relatively simple and elegant.
For the GRA method, a parameter termed the distinguishing coefficient must be prescribed.
This parameter can affect its performance; however, specific rules for assigning a value to the parameter are not available and hence, additional sensitivity analysis must be conducted.
For the DEA method, three alternatives, i.e., $\bm{x}_{11},\bm{x}_{15}$, and $\bm{x}_{18}$, are suggested, but further efforts are required to reach the final decision.
For the TOPSIS method, it leads to alternative $\bm{x}_{11}$ that is inconsistent with the consensus.
Regarding the SAW method, although it is effective in this example, limited applications have been found in the literature
because it sometimes produces results  that are not logical~\cite{velasquez2013analysis}.

\section{Conclusion}\label{sec_con}

In existing studies, a large number of MOEAs have been developed to solve MOPs. In the end a final solution must be selected out of obtained Pareto optimal solutions.
Although many MCDM approaches from the field of operations research can be adopted, they mostly require weighting coefficients prescribed by the DM and some of them lack geometric interpretations.
In the field of evolutionary computation that values geometric interpretations, few approaches to MCDM in MOPs have been developed.
In this paper, we proposed a MMD approach to MCDM in MOPs. The approach has rich geometric interpretations and avoids subjective preference inputs from the DM.
In contrast with conventional WS approaches, the MMD approach provides a systematic way to generate weighting coefficients without \emph{a priori} preference from the DM.
Our analysis showed that the approach is equivalent to knee selection described by the D\&C approach.
Simulations have been performed to illustrate the effectiveness of the proposed methodology.

\begin{table}
  \centering
  \caption{MCDM in DTLZ1 by MMD and WS Approaches}\label{tab_sim_5D}
{\scriptsize
\begin{tabular}{|c|c|c|c|}
\hline
  Solutions & $\bm{y}([\bm{x}_i])$ & MMD& WS  \\
\hline
$\bm{x}_{1}$	& $[$ 0.0074\;   0.0026\;  0.0152\;  0.1500\;  1.0080$]^T$ &1.1113&1.1833 \\
$\bm{x}_{2}$	&$[$ 0.0084\;   0.0281\;  0.0476\;  0.0830\;  0.7508$]^T$	&0.8462&0.9181	\\
$\bm{x}_{3}$	&$[$ 0.0397\;   0.0009\;  0.2390\;  0.5895\;  0.3838$]^T$	&1.1813&1.2533	\\
$\bm{x}_{4}$	&$[$ 0.0786\;   0.1104\;  0.9212\;  0.3954\;  0.3643$]^T$	&1.7981&1.8701	\\
$\bm{x}_{5}$	&$[$ 0.1045\;   0.2175\;  0.2645\;  0.8646\;  0.2316$]^T$	&1.6109&1.6828	\\
$\bm{x}_{6}$	&\textbf{$[$ 0.1075\;   0.1562\;  0.0634\;  0.0403\;  0.5492$]^T$}	&\textbf{0.8448}&\textbf{0.9167}	\\
$\bm{x}_{7}$	&$[$ 0.1081\;   0.0656\;  0.4108\;  1.0403\;  0.2550$]^T$	&1.8081&1.8800	\\
$\bm{x}_{8}$	&$[$ 0.1494\;   0.2953\;  0.1129\;  0.4294\;  0.0080$]^T$	&0.9232&0.9952	\\
$\bm{x}_{9}$	&$[$ 0.1845\;   1.0010\;  0.0744\;  0.3853\;  0.2971$]^T$	&1.8704&1.9424	\\
$\bm{x}_{10}$	&$[$ 0.1915\;   0.2743\;  1.0152\;  0.1228\;  0.3714$]^T$	&1.9035&1.9754	\\
$\bm{x}_{11}$	&$[$ 0.3801\;   0.1362\;  0.0425\;  0.7685\;  0.0800$]^T$	&1.3355&1.4075	\\
$\bm{x}_{12}$	&$[$ 0.4236\;   0.1452\;  0.5504\;  0.5501\;  0.1205$]^T$	&1.7180&1.7899	\\
$\bm{x}_{13}$	&$[$ 0.5124\;   0.7438\;  0.0866\;  0.0797\;  0.0101$]^T$	&1.3610&1.4330	\\
$\bm{x}_{14}$	&$[$ 0.6835\;   0.2687\;  0.1543\;  0.2769\;  0.1571$]^T$	&1.4688&1.5407	\\
$\bm{x}_{15}$	&$[$ 0.8185\;   0.4825\;  0.3371\;  0.2555\;  0.1091$]^T$	&1.9310&2.0029	\\
$\bm{x}_{16}$	&$[$ 1.0074\;   0.3698\;  0.1104\;  0.2089\;  0.0911$]^T$	&1.7158&1.7878	\\
\hline
\end{tabular}
}
\end{table}

\begin{figure}
  \centering
\begin{equation*}
  \begin{array}{cc}
  \hspace{-0.5cm}  \includegraphics[width=5.4cm]{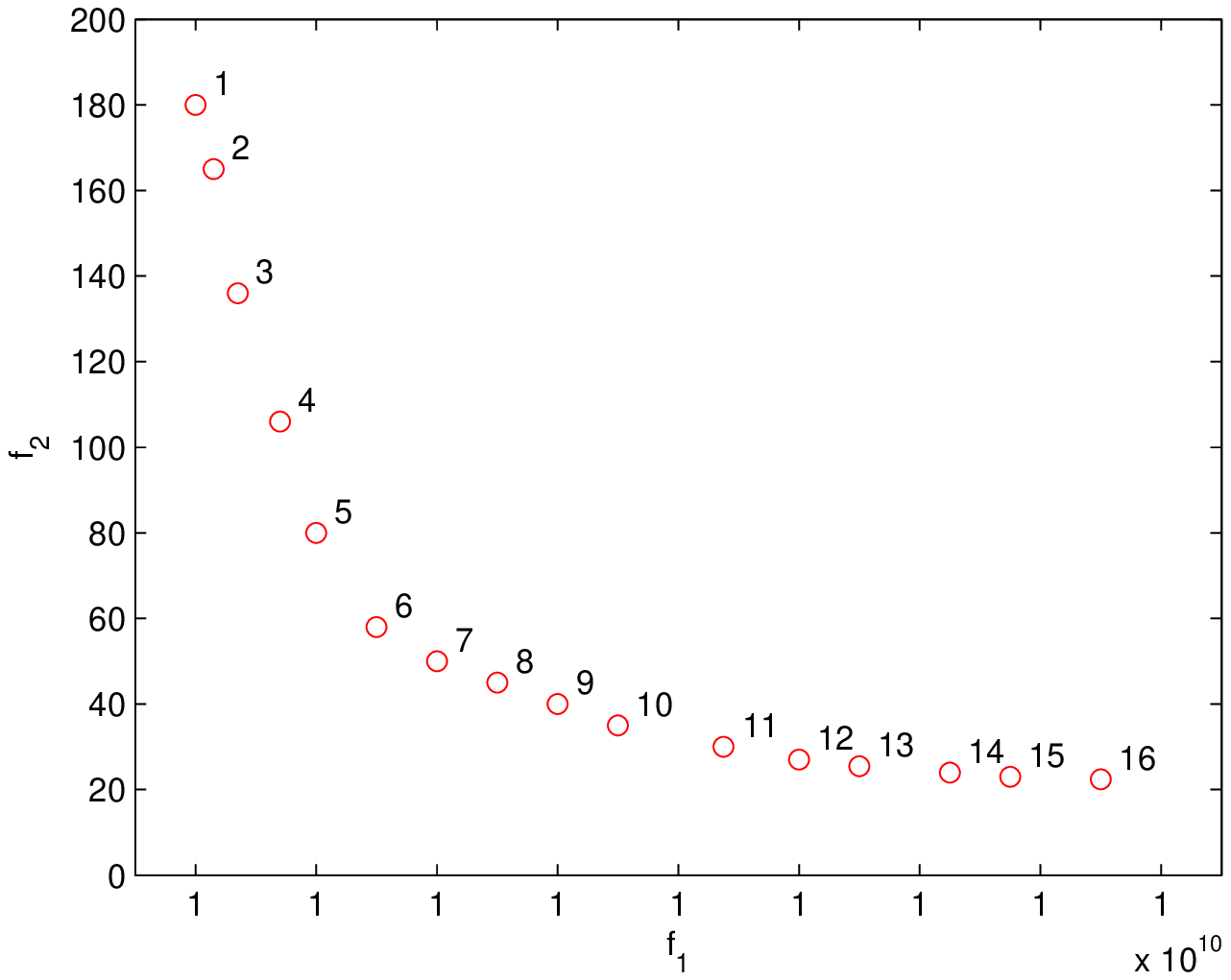}&  \includegraphics[width=3.5cm]{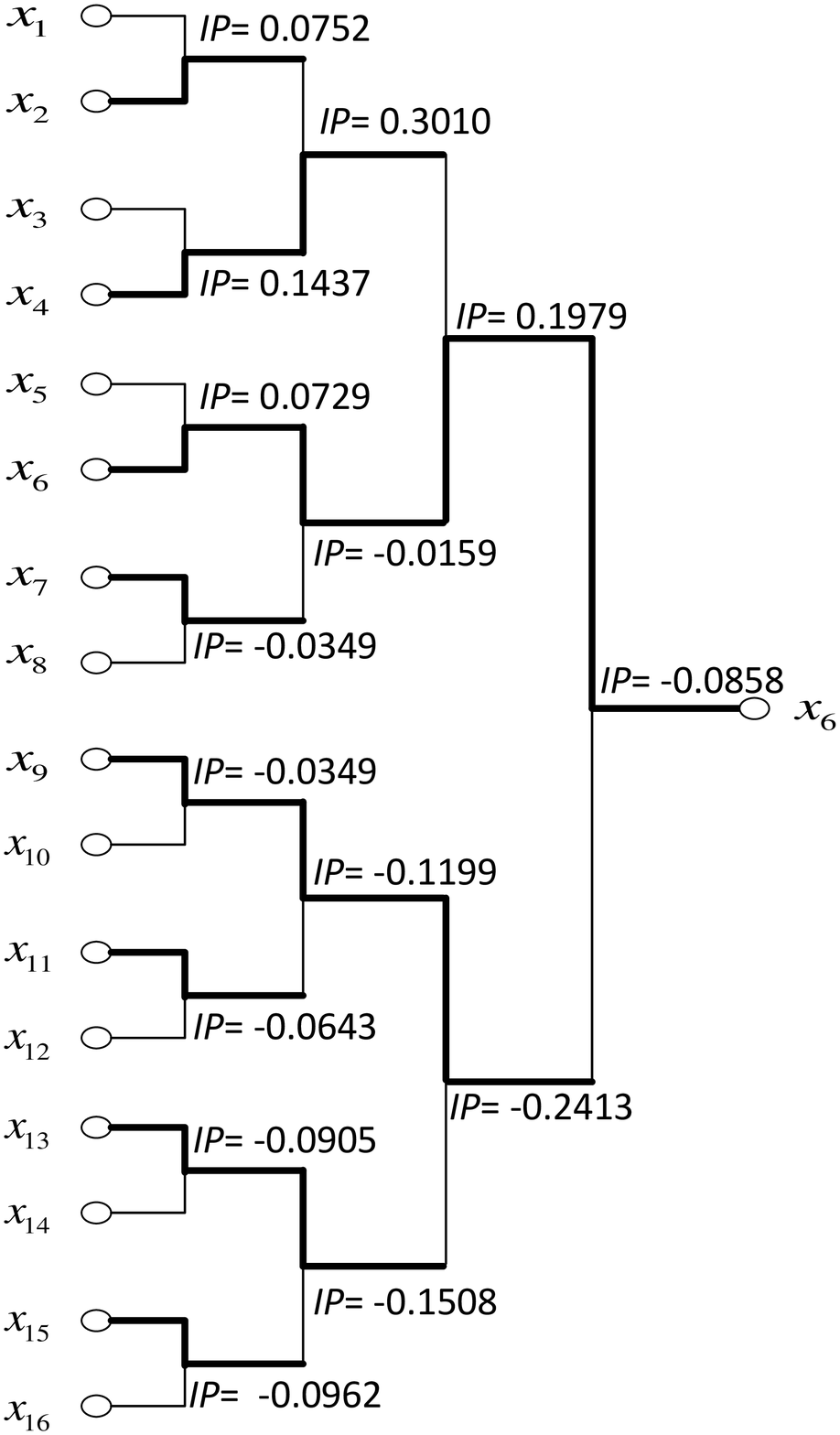} \\
   \mbox{(a)} &  \mbox{(b)}
  \end{array}
\end{equation*}
  \caption{A scenario in which all solutions yield almost the same weighted sum.
  (a) The $x$-axis is marked by the same graduation because the values of $f_1$ are much larger than the associated maximum spread  ($L_1$ is relatively small as compared to the values of $f_1$). (b) Pairwise comparisons using the D\&C approach to MCDM, where  $IP=IP(\bm{x}_i  \rightarrow \bm{x}_j )$ with $\bm{x}_i$ and $\bm{x}_j$  representing the upper and lower solutions in a pairwise comparison, respectively.
  }\label{fig_WS_fail}
\end{figure}

\begin{table}
  \centering
  \caption{MCDM in Fig.~\ref{fig_WS_fail} by MMD and WS Approaches}\label{tab_WS_fail}
\begin{tabular}{|c|c|c|}
\hline
Solutions& MMD & WS \\
\hline
$\bm{x}_{1}$    & 1&6.6667e+08  \\
\hline
$\bm{x}_{2}$    & 0.9248&6.6667e+08  \\
\hline
$\bm{x}_{3}$    & 0.7675&6.6667e+08  \\
\hline
$\bm{x}_{4}$    & 0.6238&6.6667e+08  \\
\hline
$\bm{x}_{5}$    & 0.4988&6.6667e+08  \\
\hline
$\bm{x}_{6}$    & \textbf{0.4259}&6.6667e+08  \\
\hline
$\bm{x}_{7}$    & 0.4418&6.6667e+08  \\
\hline
$\bm{x}_{8}$    & 0.4767&6.6667e+08  \\
\hline
$\bm{x}_{9}$    & 0.5117&6.6667e+08  \\
\hline
$\bm{x}_{10}$   & 0.5466&6.6667e+08  \\
\hline
$\bm{x}_{11}$   & 0.6316&6.6667e+08  \\
\hline
$\bm{x}_{12}$   & 0.6959&6.6667e+08  \\
\hline
$\bm{x}_{13}$   & 0.7530&6.6667e+08  \\
\hline
$\bm{x}_{14}$   & 0.8435&6.6667e+08  \\
\hline
$\bm{x}_{15}$   & 0.9038&6.6667e+08  \\
\hline
$\bm{x}_{16}$   & 1&6.6667e+08  \\
\hline
\end{tabular}
\end{table}

\begin{table*}
  \centering
  \caption{Comparison of Computational Time}\label{tab_sim_time}
{\tiny
\begin{tabular}{|cl|c|c|c|c|c|c|c|c|c|}
\hline
\multicolumn{2}{|c|}{\multirow{2}{*}{Problem names}} &\multicolumn{1}{c|}{ \multirow{2}{*}{$N$} }
&\multicolumn{1}{c|}{ \multirow{2}{*}{Population Size} }
&\multicolumn{1}{c|}{\multirow{2}{*}{ Num. of Simulation Runs}} &\multicolumn{2}{c|}{MMD}&\multicolumn{2}{c|}{WS}&\multicolumn{2}{c|}{D\&C}\\
\cline{6-11}
 & & & & &total time&average time&total time&average time&total time&average time\\
\hline

\multirow{6}{*}{$C_1$}&DTLZ1 &5 & 50 & 100  & 0.0082& 8.1858e-05& 0.0039& 3.9041e-05& 0.0256& 2.6559e-04\\
&DTLZ1 &5 & 50 & 1000 & 0.0482& 4.8211e-05& 0.0345& 3.4521e-05& 0.2319& 2.3191e-04\\
&DTLZ1 &5 & 50 & 3000 & 0.1114& 3.7148e-05& 0.0906& 3.0207e-05& 0.6964& 2.3213e-04\\
&\bpara{-5}{6}{180}{19}DTLZ1& 5 &50 & 5000 & 0.1746& 3.4914e-05& 0.1469& 2.9379e-05& 1.1403& 2.2806e-04\\
&DTLZ1 &5 & 50 & 8000 & 0.3005& 3.7562e-05& 0.2303& 2.8784e-05& 1.6512& 2.0640e-04\\
&DTLZ1 &5 & 50 & 10000& 0.3800& 3.8001e-05& 0.2881& 2.8807e-05& 2.3347& 2.3347e-04\\
\hline

\multirow{4}{*}{$C_2$}&DTLZ1 &5 & 25  & 3000 & 0.1009& 3.3633e-05& 0.0846& 2.8203e-05& 0.5965& 1.9883e-04\\
&DTLZ1 &5 & 50  & 3000 & 0.1114& 3.7148e-05& 0.0906& 3.0207e-05& 0.6964& 2.3213e-04\\
&\bpara{-5}{6}{180}{12}DTLZ1 &5 & 100 & 3000 & 0.1201& 4.0021e-05& 0.1036& 3.4538e-05& 1.6449& 5.4831e-04\\
&DTLZ1 &5 & 200 & 3000 & 0.1400& 4.6652e-05& 0.1148& 3.8281e-05& 1.9489& 6.4962e-04\\
\hline

\multirow{16}{*}{$C_3$}&DTLZ1 &5& 50  & 3000 & 0.1114& 3.7148e-05& 0.0906& 3.0207e-05& 0.6964& 2.3213e-04\\
&DTLZ2 &5& 50  & 3000 & 0.1258& 4.1927e-05& 0.0967& 3.2295e-05& 0.9329& 3.1097e-04\\
&MOP1  &2& 50  & 3000 & 0.2593& 8.6417e-05& 0.0930& 3.1001e-05& 0.9492& 3.1639e-04\\
&MOP2  &2& 50  & 3000 & 0.1163& 3.8777e-05& 0.0842& 2.8083e-05& 0.5967& 1.9890e-04\\
&MOP3  &2& 50  & 3000 & 0.1122& 3.7411e-05& 0.0835& 2.7841e-05& 0.5923& 1.9744e-04\\
&MOP4  &2& 50  & 3000 & 0.1212& 4.0411e-05& 0.0885& 2.9494e-05& 0.9326& 3.1088e-04\\
&MOP5  &3& 50  & 3000 & 0.1230& 4.1010e-05& 0.0907& 3.0228e-05& 0.9384& 3.1279e-04\\
&MOP6  &2& 50  & 3000 & 0.1177& 3.9237e-05& 0.0857& 2.8555e-05& 0.8660& 2.8866e-04\\
&\bpara{-5}{6}{180}{50}MOP7  & 3  &50& 3000 & 0.1272& 4.2416e-05& 0.0915& 3.0503e-05& 0.9529& 3.1764e-04\\
&MOP-C1 Binh     &2& 50  & 3000 & 0.1222& 4.0738e-05& 0.0875& 2.9157e-05& 0.9222& 3.0741e-04\\
&MOP-C1 Osyczka  &2& 50  & 3000 & 0.1209& 4.0304e-05& 0.0880& 2.9339e-05& 0.9032& 3.0105e-04\\
&MOP-C1 Viennet  &3& 50  & 3000 & 0.1242& 4.1414e-05& 0.0929& 3.0951e-05& 0.9237& 3.0789e-04\\
&MOP-C1 Tanaka   &2& 50  & 3000 & 0.1216& 4.0531e-05& 0.0873& 2.9100e-05& 0.9115& 3.0383e-04\\
&ZDT1  &2& 50  & 3000 & 0.1433& 4.7776e-05& 0.1009& 3.3644e-05& 0.7088& 2.3626e-04\\
&ZDT2  &2& 50  & 3000 & 0.1155& 3.8508e-05& 0.0825& 2.7515e-05& 0.5166& 1.7219e-04\\
&ZDT3  &2& 50  & 3000 & 0.1113& 3.7102e-05& 0.0830& 2.7669e-05& 0.5619& 1.8729e-04\\
\hline
\end{tabular}
}
\end{table*}

\begin{figure*}
  \centering
\begin{equation*}
\begin{array}{ccc}
\includegraphics[width=5cm]{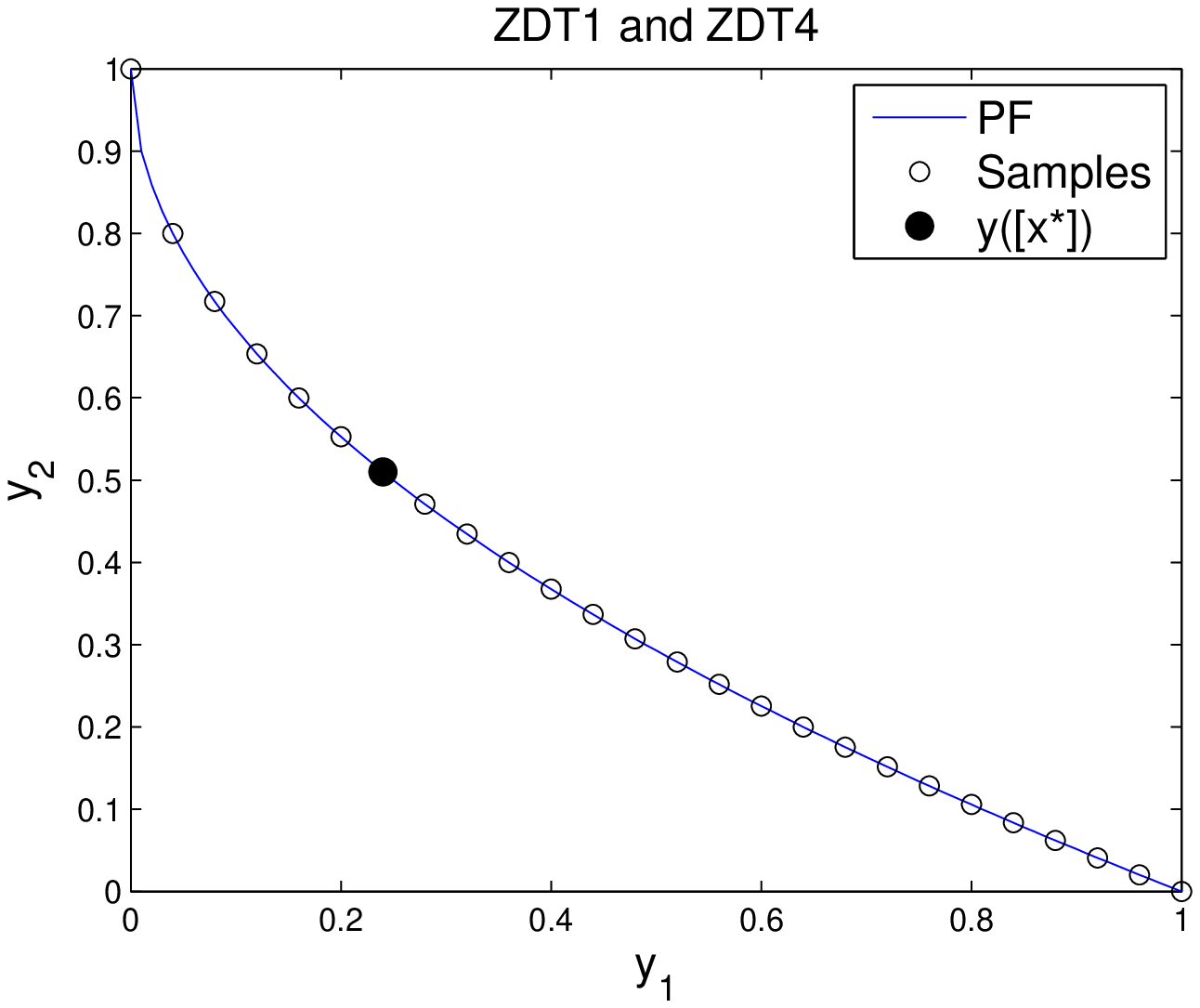} &  \includegraphics[width=5cm]{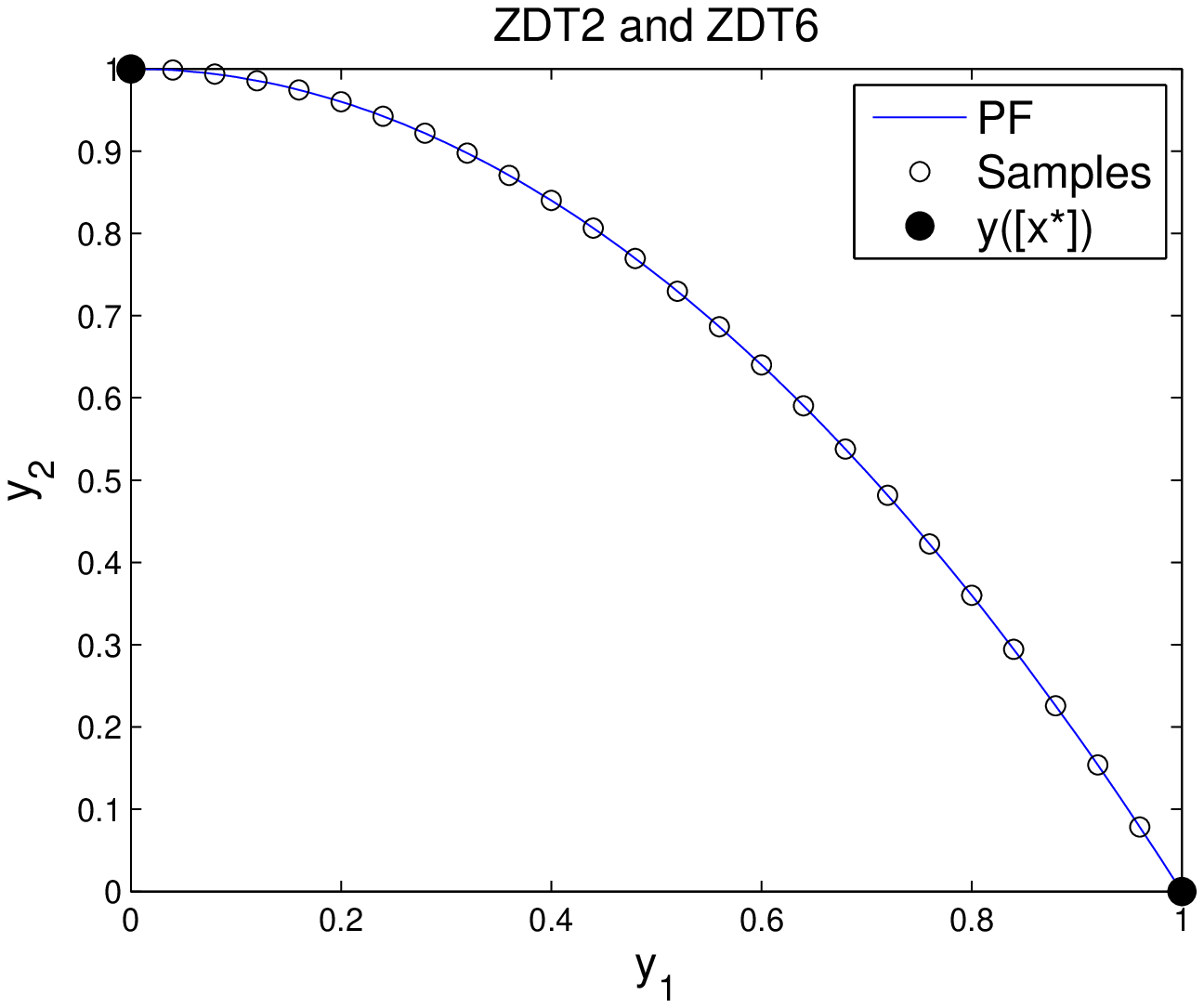}  & \includegraphics[width=5cm]{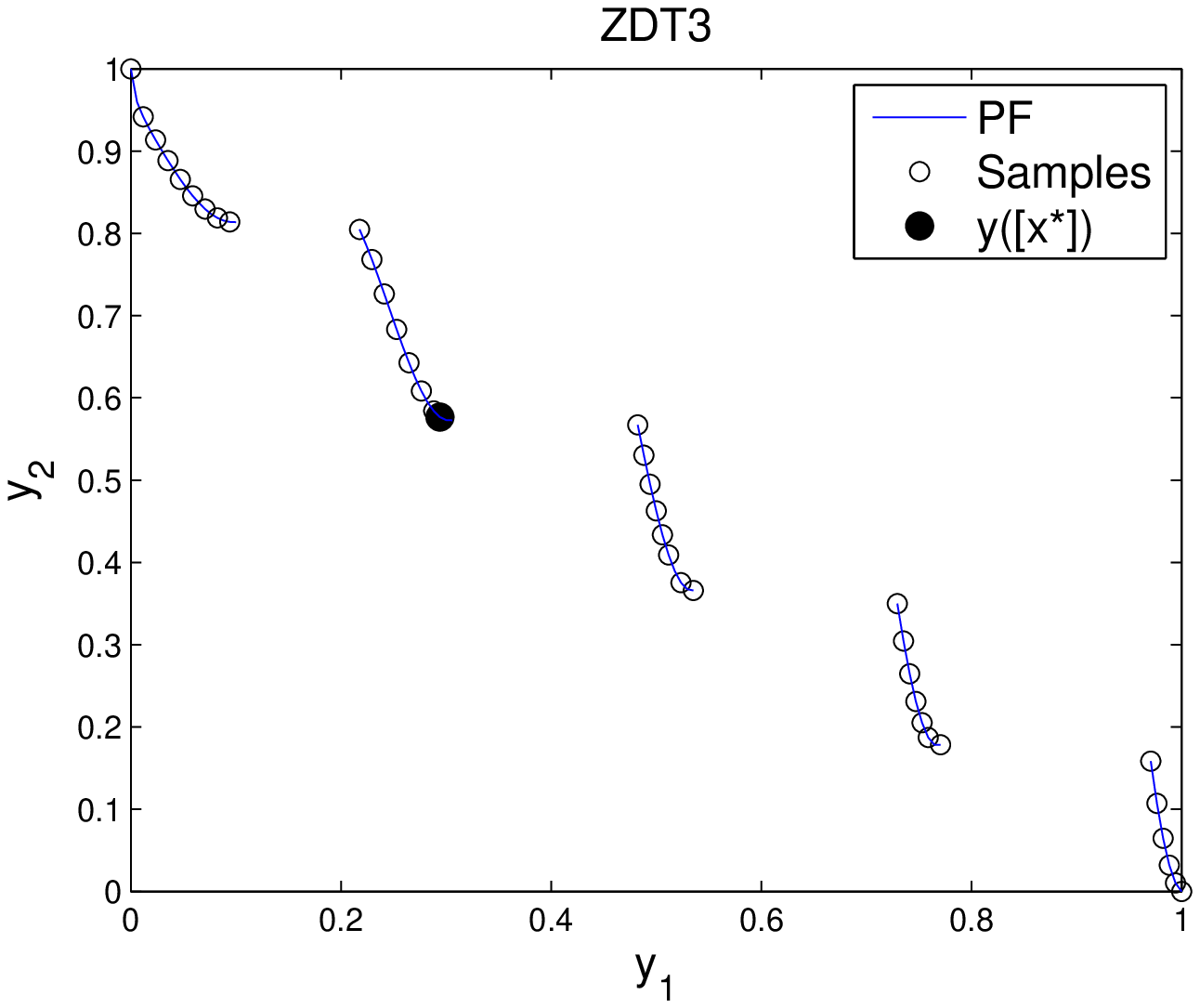} \\
  \mbox{(a)} & \mbox{(b)} & \mbox{(c)}  \\
  \includegraphics[width=5cm]{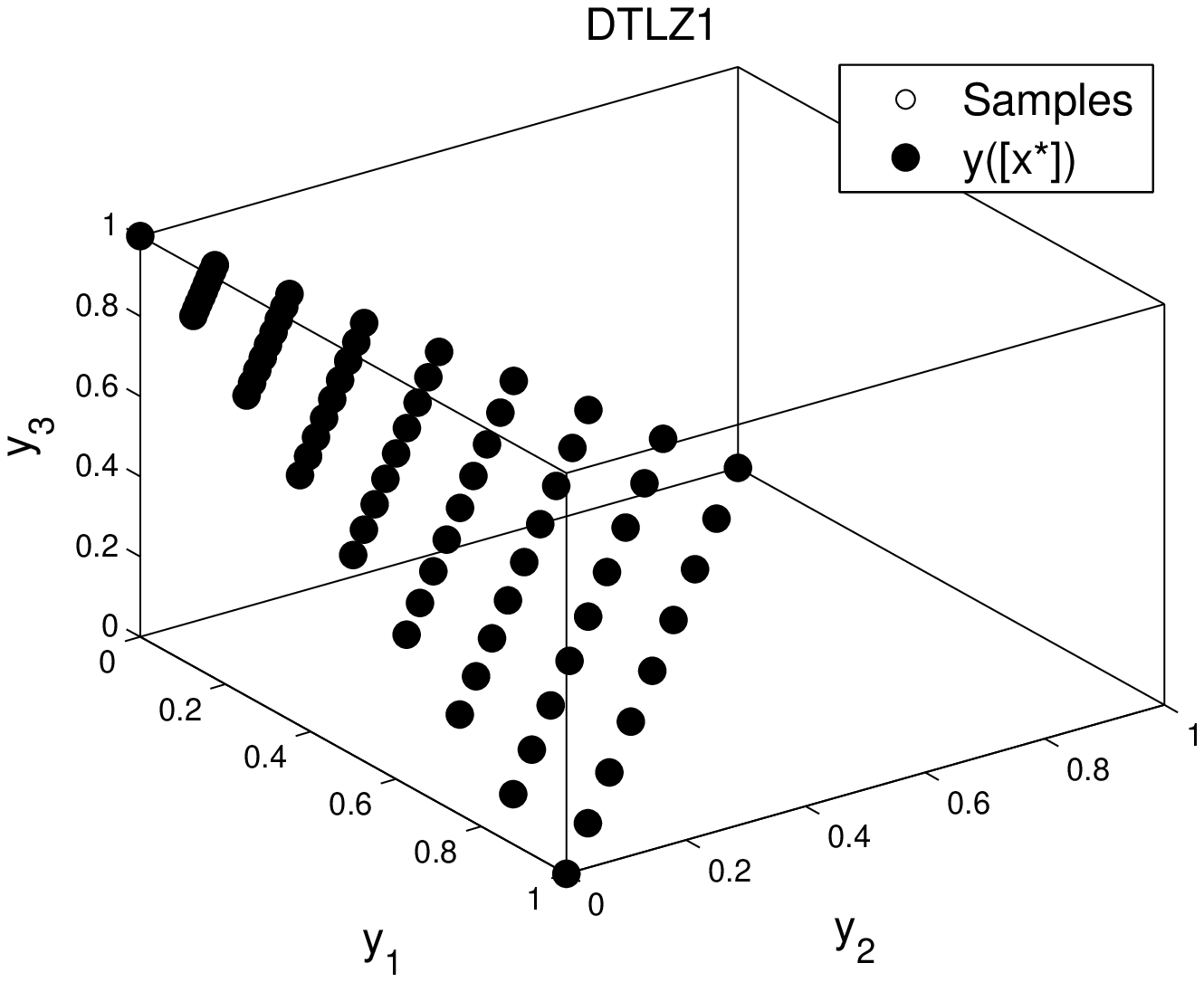}  & \includegraphics[width=5cm]{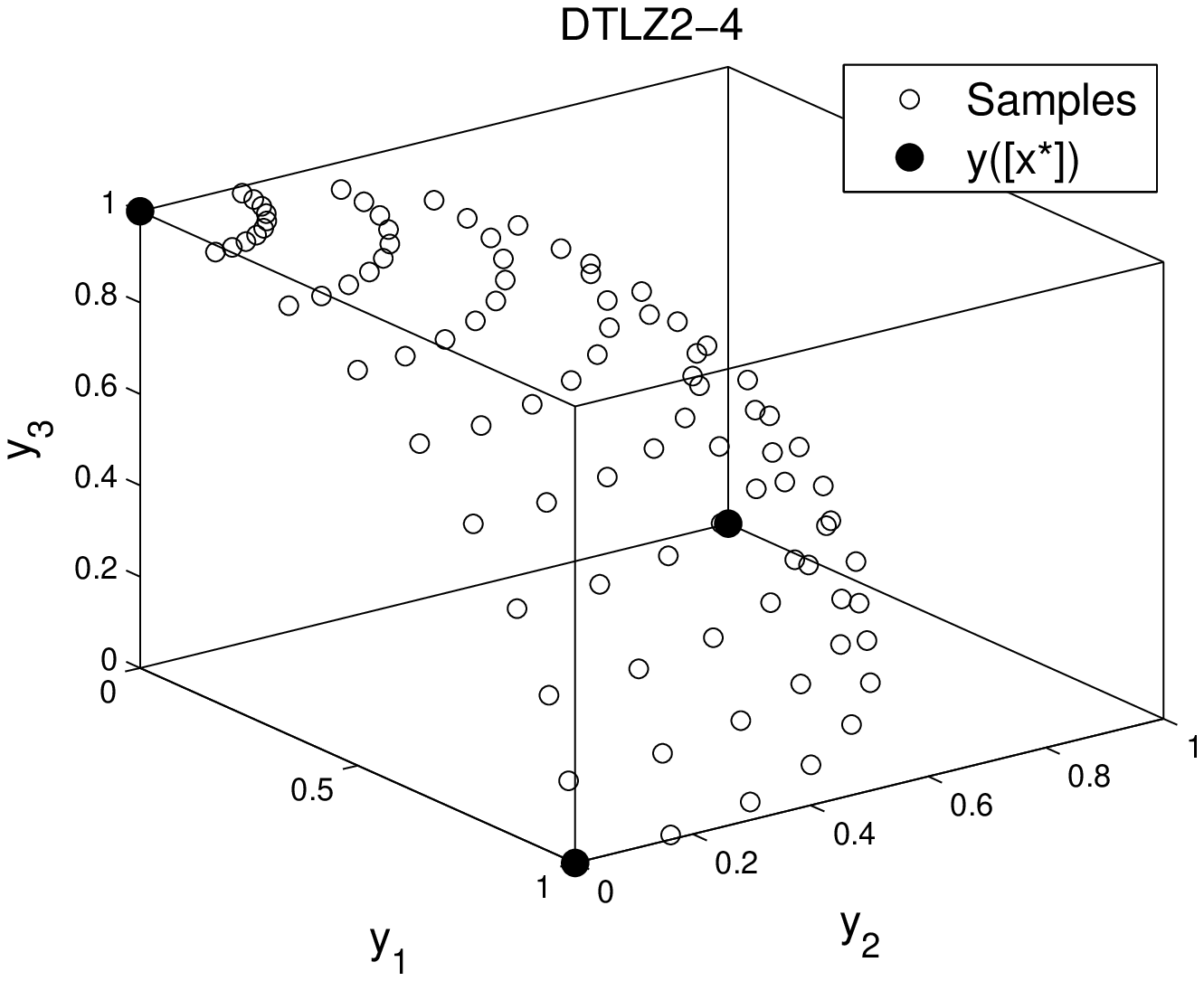} &   \includegraphics[width=5cm]{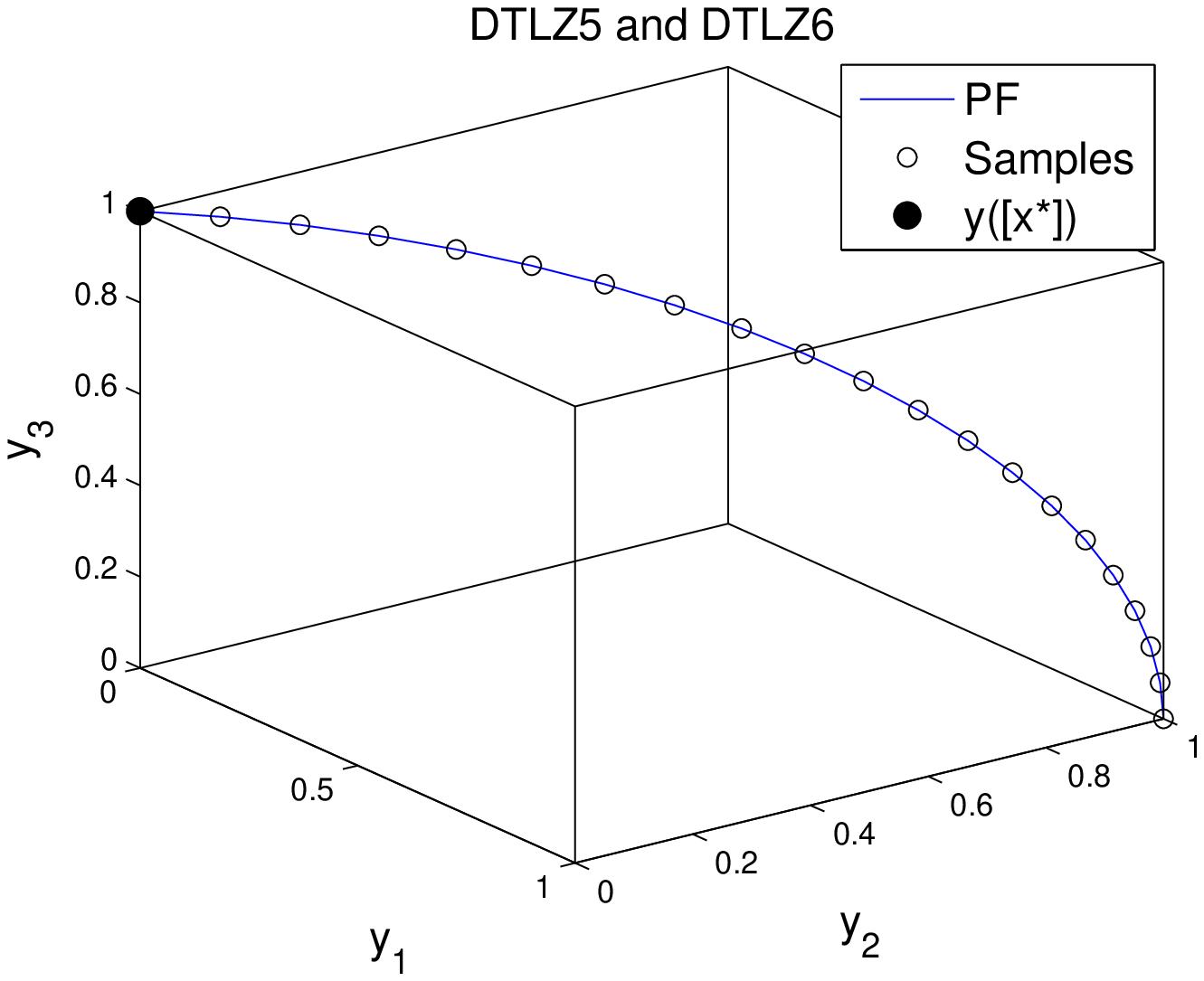} \\
  \mbox{(d)} & \mbox{(e)} & \mbox{(f)}  \\
 \includegraphics[width=5cm]{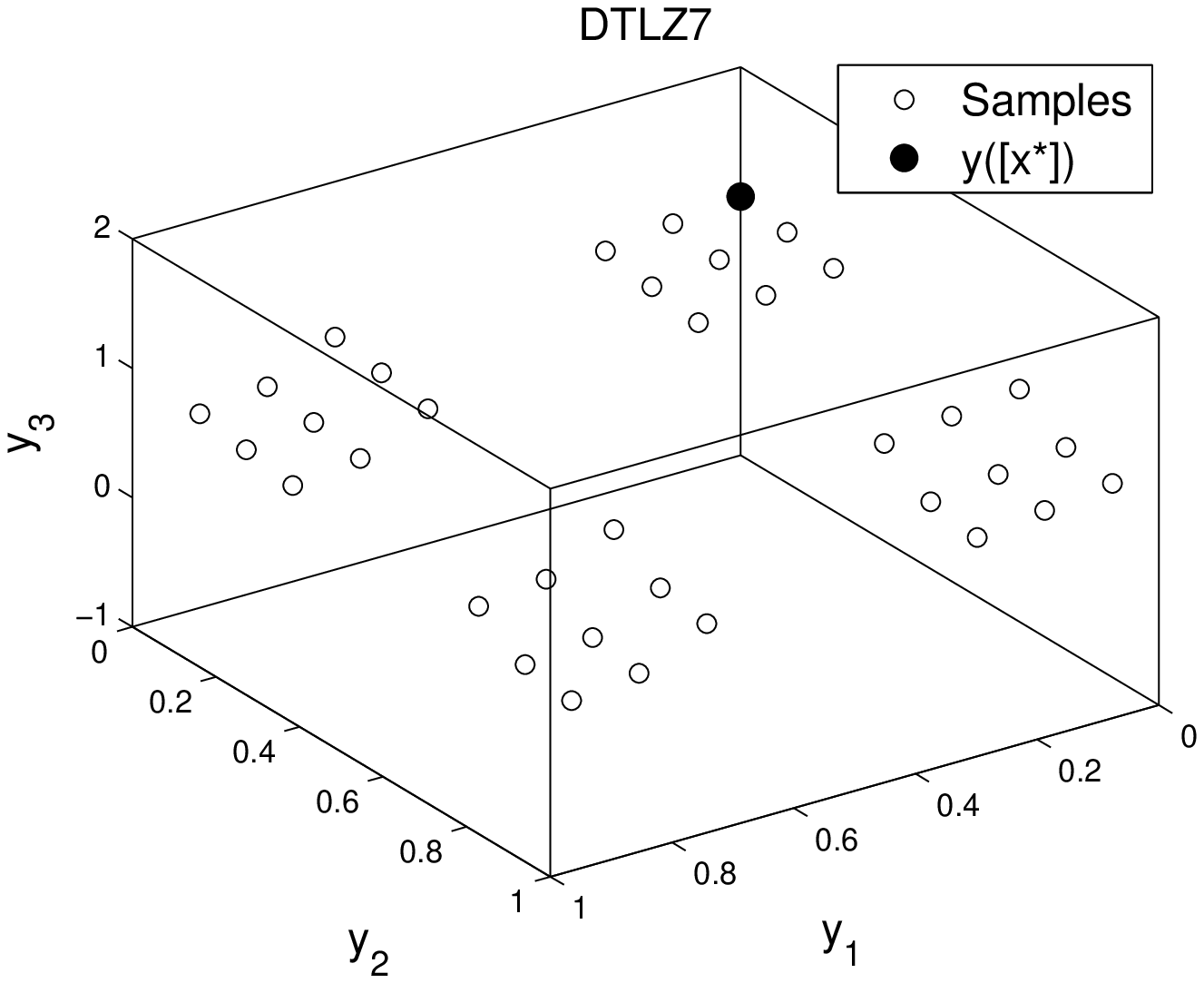}   &  \includegraphics[width=5cm]{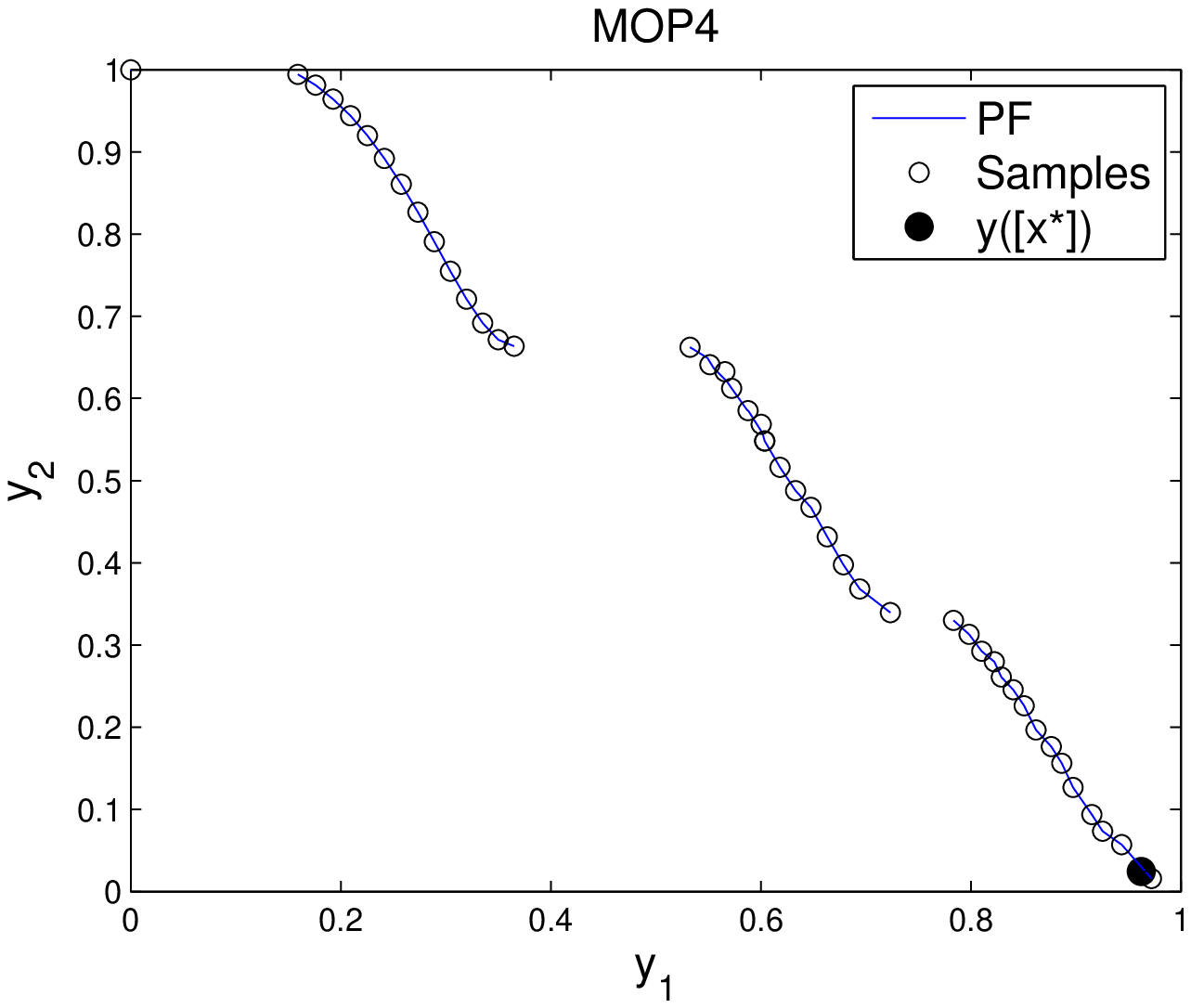}  &  \includegraphics[width=5cm]{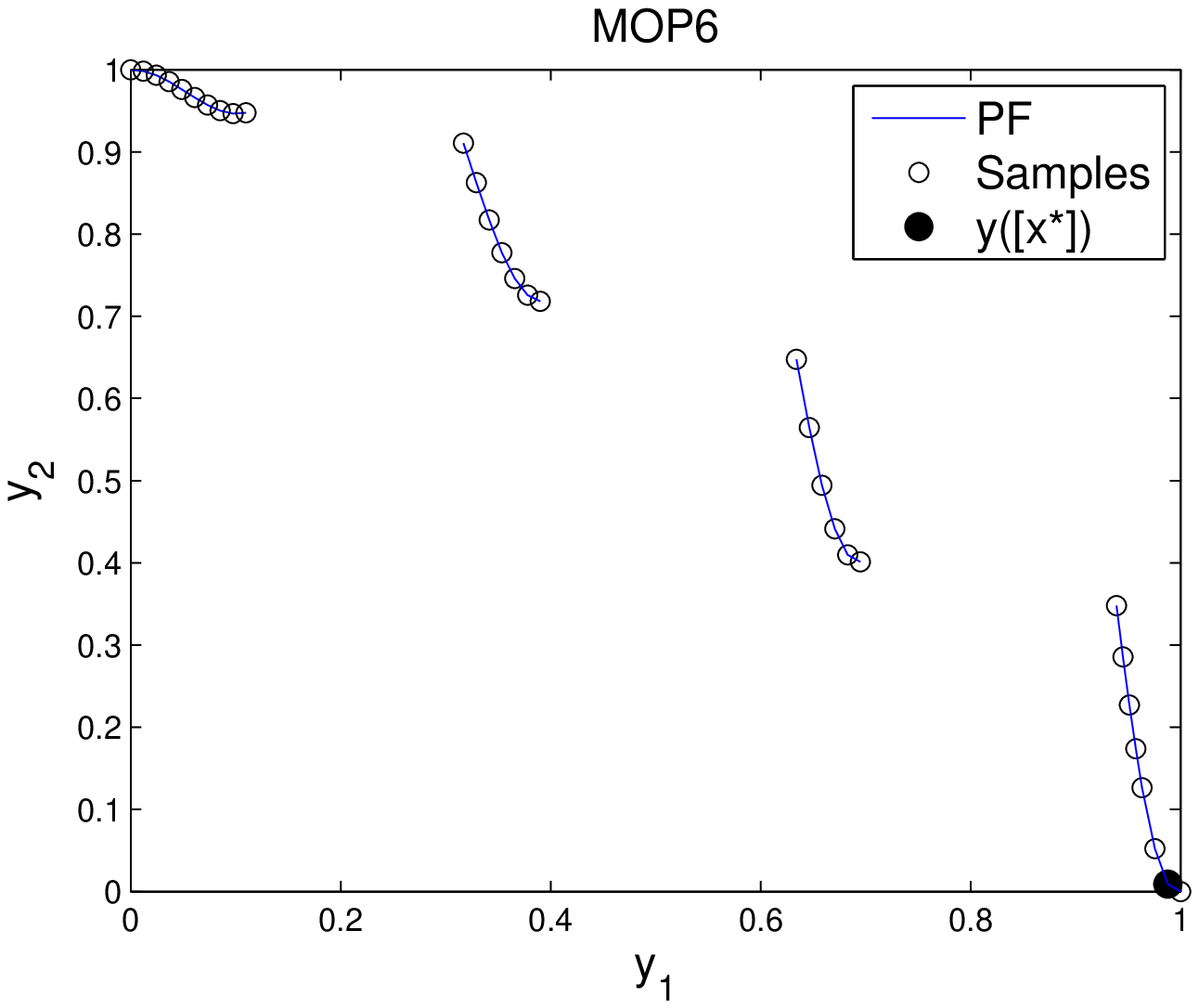}  \\
  \mbox{(g)} & \mbox{(h)} & \mbox{(i)} \\
\end{array}
\end{equation*}
  \caption{MCDM in ZDT and DTLZ test suites, MOP4, and MOP6 using the MMD approach. In a 2-D space, the PF in (a) has the shape of a convex curve, and
  the PF in (b) has the shape of a concave curve.
  In a 3-D space, the PF in (d) is on a plane, and the PF in (e)  has the shape of a concave surface.
  Discontinuous PFs appear in (c), (g), (h), and (i). The selection process can be readily visualized because of the front shapes except for the PFs in (f) and (g). }\label{fig_ZDT_DTLZ}
\end{figure*}

\begin{table}
  \centering
  \caption{MCDM for Facility Layout} \label{tab_RealData}
\begin{tabular}{cccccc}
\hline
Ranking &GRA  &DEA &TOPSIS &SAW &MMD\\
&  \cite{kuo2008use} & \cite{yang03hierarchical} & \cite{yang2007multiple}  &  \cite{kuo2008use} & \\
\hline  \hline

1& $\bm{x}_{15}$ & $\bm{x}_{11}$ $\bm{x}_{15}$ $\bm{x}_{18}$& $\bm{x}_{11}$& $\bm{x}_{15}$& $\bm{x}_{15}$\\

\hline
2& $\bm{x}_{17}$& & $\bm{x}_{15}$& $\bm{x}_{17}$& $\bm{x}_{11}$\\
\hline
3& $\bm{x}_{11}$& & $\bm{x}_{17}$& $\bm{x}_{11}$& $\bm{x}_{17}$\\
\hline
4& $\bm{x}_{18}$& $\bm{x}_{2}$& $\bm{x}_{16}$& $\bm{x}_{18}$& $\bm{x}_{18}$\\
\hline
5& $\bm{x}_{9}$& $\bm{x}_{16}$& $\bm{x}_{10}$& $\bm{x}_{9}$& $\bm{x}_{9}$\\
\hline
6& $\bm{x}_{16}$& $\bm{x}_{6}$& $\bm{x}_{9}$& $\bm{x}_{16}$& $\bm{x}_{16}$\\
\hline
7& $\bm{x}_{8}$& $\bm{x}_{8}$& $\bm{x}_{18}$& $\bm{x}_{10}$& $\bm{x}_{2}$\\
\hline
8& $\bm{x}_{2}$& $\bm{x}_{9}$& $\bm{x}_{2}$& $\bm{x}_{2}$& $\bm{x}_{10}$\\
\hline
9& $\bm{x}_{10}$& $\bm{x}_{17}$& $\bm{x}_{6}$& $\bm{x}_{8}$& $\bm{x}_{8}$\\
\hline
10& $\bm{x}_{1}$& $\bm{x}_{1}$& $\bm{x}_{13}$& $\bm{x}_{14}$& $\bm{x}_{14}$\\
\hline\hline
\end{tabular}
\end{table}

\section*{Acknowledgment}
The authors would like to thank Dr. Zhenan He for generating the Pareto fronts in Fig.~\ref{fig_ZDT_DTLZ}.


\begin{thebibliography}{10}
\providecommand{\url}[1]{#1}
\csname url@samestyle\endcsname
\providecommand{\newblock}{\relax}
\providecommand{\bibinfo}[2]{#2}
\providecommand{\BIBentrySTDinterwordspacing}{\spaceskip=0pt\relax}
\providecommand{\BIBentryALTinterwordstretchfactor}{4}
\providecommand{\BIBentryALTinterwordspacing}{\spaceskip=\fontdimen2\font plus
\BIBentryALTinterwordstretchfactor\fontdimen3\font minus
  \fontdimen4\font\relax}
\providecommand{\BIBforeignlanguage}[2]{{%
\expandafter\ifx\csname l@#1\endcsname\relax
\typeout{** WARNING: IEEEtran.bst: No hyphenation pattern has been}%
\typeout{** loaded for the language `#1'. Using the pattern for}%
\typeout{** the default language instead.}%
\else
\language=\csname l@#1\endcsname
\fi
#2}}
\providecommand{\BIBdecl}{\relax}
\BIBdecl

\bibitem{Hamming}
M.~Izadikhah, ``Using the {H}amming distance to extend {TOPSIS} in a fuzzy
  environment,'' \emph{Journal of Computational and Applied Mathematics}, vol.
  231, no.~1, pp. 200--207, 2009.

\bibitem{path_follow}
K.-B. Lee and J.-H. Kim, ``Multiobjective particle swarm optimization with
  preference-based sort and its application to path following footstep
  optimization for humanoid robots,'' \emph{{IEEE} Trans. Evol. Comput.},
  vol.~17, no.~6, pp. 755--766, Dec. 2013.

\bibitem{factor_effi}
Y.~Kuo, T.~Yang, and G.-W. Huang, ``The use of grey relational analysis in
  solving multiple attribute decision-making problems,'' \emph{Computer \&
  Industrial Engineering}, vol.~55, pp. 80--93, 2008.

\bibitem{ORs}
H.~A. Eiselt and C.-L. Sandblom, \emph{Operations Research: A Model-based
  Approach}.\hskip 1em plus 0.5em minus 0.4em\relax Berlin, Heidelberg:
  Springer, 2012.

\bibitem{Knowles99PAES}
J.~D. Knowles and D.~W. Corne, ``The {P}areto archived evolution strategy: a
  new baseline algorithm for {P}areto multiobjective optimization,'' in
  \emph{Proc. IEEE Conf. Evolutionary Computation}, Washington, DC, Jul. 1999,
  pp. 98--105.

\bibitem{Corne00PESA}
D.~W. Corne, J.~D. Knowles, and M.~J. Oates, ``The {P}areto envelope-based
  selection algorithm for multiobjective optimization,'' in \emph{Proc. Int.
  Conf. Parallel Problem Solving From Nature}, Paris, France, Sep. 2000, pp.
  839--848.

\bibitem{Deb02NSGA-II}
K.~Deb, S.~Agrawal, A.~Pratap, and T.~Meyarivan, ``A fast and elitist
  multiobjective genetic algorithm: {NSGA-II},'' \emph{{IEEE} Trans. Evol.
  Comput.}, vol.~6, no.~2, pp. 182--197, Apr. 2002.

\bibitem{Zitzler02SPEA2}
E.~Zitzler, M.~Laumanns, and L.~Thiele, ``{SPEA2}: improving the strength
  pareto evolutionary algorithm,'' in \emph{Proc. Evol. Methods Design
  Optimization Control Applicat. Ind. Problems}, Athens, Greece, Apr. 2002, pp.
  95--100.

\bibitem{Zhang07MOEA}
Q.~Zhang and H.~Li, ``{MOEA/D}: a multiobjective evolutionary algorithm based
  on decomposition,'' \emph{{IEEE} Trans. Evol. Comput.}, vol.~11, no.~6, pp.
  712--731, Dec. 2007.

\bibitem{three_category}
D.~{Van Veldhuizen} and G.~B. Lamont, ``Multiobjective evolutionary algorithms:
  analyzing the state-of-the-art,'' \emph{Evolutionary Computation Journal},
  vol.~8, no.~2, pp. 125--147, 2000.

\bibitem{decom_series}
C.~Audet, G.~Savard, and W.~Zghal, ``Multiobjective optimization through a
  series of single-objective formulations,'' \emph{SIAM J. Optim.}, vol.~19,
  no.~1, pp. 188--210, 2008.

\bibitem{Rachmawati06}
L.~Rachmawati and D.~Srinivasan, ``A multi-objective genetic algorithm with
  controllable convergence on knee regions,'' in \emph{Proc. IEEE Congress on
  Evolutionary Computation}, Vancouver, British Columbia, Canada, Jul. 2006,
  pp. 1916--1923.

\bibitem{pre_based}
J.-H. Kim, J.-H. Han, Y.-H. Kim, S.-H. Choi, and E.-S. Kim, ``Preference-based
  solution selection algorithm for evolutionary multiobjective optimization,''
  \emph{{IEEE} Trans. Evol. Comput.}, vol.~16, no.~1, pp. 20--34, Feb. 2012.

\bibitem{branke2004finding}
J.~Branke, K.~Deb, H.~Dierolf, and M.~Osswald, ``Finding knees in
  multi-objective optimization,'' in \emph{Proc. Int. Conf. Parallel Problem
  Solving from Nature (PPSN)}, Birmingham, UK, Sep. 2004, pp. 722--731.

\bibitem{Cvet2002pre}
D.~Cvetkovic and I.~C. Parmee, ``Preferences and their application in
  evolutionary multiobjective optimization,'' \emph{{IEEE} Trans. Evol.
  Comput.}, vol.~6, no.~1, pp. 42--57, Feb. 2002.

\bibitem{Rach09knee}
L.~Rachmawati and D.~Srinivasan, ``Multiobjective evolutionary algorithm with
  controllable focus on the knees of the pareto front,'' \emph{{IEEE} Trans.
  Evol. Comput.}, vol.~13, no.~4, pp. 810--824, Aug. 2009.

\bibitem{knee_progress}
X.~Zhang, Y.~Tian, and Y.~Jin, ``A knee point-driven evolutionary algorithm for
  many-objective optimization,'' \emph{{IEEE} Trans. Evol. Comput.}, vol.~19,
  no.~6, pp. 761--776, Dec. 2015.

\bibitem{marler2004survey}
R.~T. Marler and J.~S. Arora, ``Survey of multi-objective optimization methods
  for engineering,'' \emph{Structural and multidisciplinary optimization},
  vol.~26, no.~6, pp. 369--395, 2004.

\bibitem{das1997closer}
I.~Das and J.~E. Dennis, ``A closer look at drawbacks of minimizing weighted
  sums of objectives for {P}areto set generation in multicriteria optimization
  problems,'' \emph{Structural Optimization}, vol.~14, no.~1, pp. 63--69, 1997.

\bibitem{coit2004system}
D.~W. Coit, T.~Jin, and N.~Wattanapongsakorn, ``System optimization with
  component reliability estimation uncertainty: a multi-criteria approach,''
  \emph{IEEE Trans. Reliability}, vol.~53, no.~3, pp. 369--380, Sep. 2004.

\bibitem{visual_Tusar}
T.~Tusar and B.~Filipic, ``Visualization of {P}areto front approximations in
  evolutionary multiobjective optimization: A critical review and the
  prosection method,'' \emph{{IEEE} Trans. Evol. Comput.}, vol.~19, no.~2, pp.
  225--245, Apr. 2015.

\bibitem{geo_im}
S.~B. Gee, K.~C. Tan, V.~A. Shim, and N.~Pal, ``Online diversity assessment in
  evolutionary multiobjective optimization: A geometrical perspective,''
  \emph{{IEEE} Trans. Evol. Comput.}, vol.~19, no.~4, pp. 542--559, Aug. 2015.

\bibitem{metric}
C.-K. Goh and K.~C. Tan, \emph{Evolutionary Multi-objective Optimization in
  Uncertain Environments: Issues and Algorithms}.\hskip 1em plus 0.5em minus
  0.4em\relax Berlin, Heidelberg: Springer Berlin Heidelberg, 2009.

\bibitem{ParCoor}
A.~Inselberg, \emph{Visual Multidimensional Geometry and Its
  Applications}.\hskip 1em plus 0.5em minus 0.4em\relax New York, NY: Springer,
  2009.

\bibitem{para_cor1}
K.~Deb and H.~Jain, ``An evolutionary many-objective optimization algorithm
  using reference-point-based nondominated sorting approach, part {I}: solving
  problems with box constraints,'' \emph{{IEEE} Trans. Evol. Comput.}, vol.~18,
  no.~4, pp. 577--601, Aug. 2014.

\bibitem{para_cor2}
H.~Jain and K.~Deb, ``An evolutionary many-objective optimization algorithm
  using reference-point based nondominated sorting approach, part {II}:
  handling constraints and extending to an adaptive approach,'' \emph{{IEEE}
  Trans. Evol. Comput.}, vol.~18, no.~4, pp. 602--622, Aug. 2014.

\bibitem{Pryke15heatmap}
A.~Pryke, S.~Mostaghim, and A.~Nazemi, ``Heatmap visualisation of population
  based multi objective algorithms,'' in \emph{Proc. Int. conf. Evolutionary
  Multi-Criterion Optimization}, Matsushima, Japan, Jan. 2007, pp. 361--375.

\bibitem{Barton07SammonMapping}
J.~J. Valdes and A.~J. Barton, ``Visualizing high dimensional objective spaces
  for multi-objective optimization: a virtual reality approach,'' in
  \emph{Proc. IEEE Conf. Evolutionary Computation}, Singapore, Sep. 2007, pp.
  4199--4206.

\bibitem{Hoffman97RadialCoordinateVisualization}
P.~Hoffman, G.~Grinstein, K.~Marx, I.~Grosse, and E.~Stanley, ``{DNA} visual
  and analytic data mining,'' in \emph{Proc. IEEE Visualization Conf.},
  Phoenix, AZ, Oct. 1997, pp. 437--441.

\bibitem{He15polar}
Z.~He and G.~G. Yen, ``Visualization and performance metric in many-objective
  optimization,'' \emph{IEEE Trans. Evolutionary Computation}, 2015, early
  access.

\bibitem{Kohonen01selforganizingmap}
T.~Kohonen, ``The self-organizing map,'' \emph{Proceedings of the IEEE},
  vol.~78, no.~9, pp. 1464--1480, Sep. 1990.

\bibitem{SelfOrganizingMaps}
------, \emph{Self-Organizing Maps}.\hskip 1em plus 0.5em minus 0.4em\relax
  Berlin Heidelberg, Germany: Springer, 2001.

\bibitem{Tenenbaum00isomap}
J.~B. Tenenbaum, V.~Silva, and J.~C. Langford, ``A global geometric framework
  for nonlinear dimensionality reduction,'' \emph{Science}, vol. 290, no. 5500,
  pp. 2319--2323, Dec. 2000.

\bibitem{boyd_linear}
S.~P. Boyd and C.~H. Barratt, \emph{Linear Controller Design: Limits of
  Performance}.\hskip 1em plus 0.5em minus 0.4em\relax Englewood Cliffs, N.J.:
  Prentice Hall, 1991.

\bibitem{Mie99knee}
K.~Miettinen, \emph{Nonlinear Multiobjective Optimization}.\hskip 1em plus
  0.5em minus 0.4em\relax Boston: Kluwer Academic Publishers, 1999.

\bibitem{Matt04filter}
C.~A. Mattson, A.~A. Mullur, and A.~Messac, ``Smart {P}areto filter: obtaining
  a minimal representation of multiobjective design space,'' \emph{Eng.
  Optimization}, vol.~36, no.~6, pp. 271--740, 2004.

\bibitem{deb2003multi}
K.~Deb, ``Multi-objective evolutionary algorithms: introducing bias among
  {P}areto-optimal solutions,'' in \emph{Advances in Evolutionary
  Computing}.\hskip 1em plus 0.5em minus 0.4em\relax London, U.K.:
  Springer-Verlag, 2003, pp. 263--292.

\bibitem{IET_CTA_14}
W.-Y. Chiu, ``Multiobjective controller design by solving a multiobjective
  matrix inequality problem,'' \emph{IET Control Theory Appl.}, vol.~8, no.~16,
  pp. 1656--1665, Nov. 2014.

\bibitem{das1999characterizing}
I.~Das, ``On characterizing the ``knee'' of the {P}areto curve based on
  normal-boundary intersection,'' \emph{Structural Optimization}, vol.~18, no.
  2-3, pp. 107--115, Oct. 1999.

\bibitem{TSG_15}
W.-Y. Chiu, H.~Sun, and H.~Poor, ``A multiobjective approach to multimicrogrid
  system design,'' \emph{IEEE Trans. Smart Grid}, vol.~6, no.~5, pp.
  2263--2272, Sep. 2015.

\bibitem{MOEA_bk1}
C.~A. {Coello Coello}, D.~A. {Van Veldhuizen}, and G.~B. Lamont,
  \emph{Evolutionary Algorithms for Solving Multi-objective Problems}.\hskip
  1em plus 0.5em minus 0.4em\relax New York: Kluwer Academic, 2002.

\bibitem{MOEA_bk2}
K.~Deb, \emph{Multi-Objective Optimization Using Evolutionary
  Algorithms}.\hskip 1em plus 0.5em minus 0.4em\relax New York: Wiley, 2001.

\bibitem{set_theory}
C.~C. Pinter, \emph{Set Theory}.\hskip 1em plus 0.5em minus 0.4em\relax
  Reading, Mass.: Addison-Wesley Pub. Co., 1971.

\bibitem{algebra_set}
J.~B. Fraleigh, \emph{A First Course in Abstract Algebra}.\hskip 1em plus 0.5em
  minus 0.4em\relax Boston: Addison-Wesley, 2003.

\bibitem{pohekar2004CP}
S.~D. Pohekar and M.~Ramachandran, ``Application of multi-criteria decision
  making to sustainable energy planning--a review,'' \emph{Renewable and
  sustainable energy reviews}, vol.~8, no.~4, pp. 365--381, 2004.

\bibitem{Zeleny82CP}
M.~Zeleny, \emph{Multiple Criteria Decision Making}.\hskip 1em plus 0.5em minus
  0.4em\relax New York: McGraw-Hill, 1982.

\bibitem{kuo2008use}
Y.~Kuo, T.~Yang, and G.-W. Huang, ``The use of grey relational analysis in
  solving multiple attribute decision-making problems,'' \emph{Computers \&
  Industrial Engineering}, vol.~55, no.~1, pp. 80--93, 2008.

\bibitem{yang03hierarchical}
T.~Yang and C.~Kuo, ``A hierarchical {AHP/DEA} methodology for the facilities
  layout design problem,'' \emph{European Journal of Operational Research},
  vol. 147, no.~1, pp. 128--136, 2003.

\bibitem{yang2007multiple}
T.~Yang and C.-C. Hung, ``Multiple-attribute decision making methods for plant
  layout design problem,'' \emph{Robotics and Computer-Integrated
  Manufacturing}, vol.~23, no.~1, pp. 126--137, 2007.

\bibitem{velasquez2013analysis}
M.~Velasquez and P.~T. Hester, ``An analysis of multi-criteria decision making
  methods,'' \emph{International Journal of Operations Research}, vol.~10,
  no.~2, pp. 56--66, 2013.

\end{thebibliography}
\end{document}